\documentclass[12pt,a4paper]{article}

\usepackage{amssymb}
\usepackage{xypic}

\input xy

\xyoption{all}

\newcommand{\stl}{\vspace{3mm}}

\newcommand{\gr}{\mbox{gr}}

\newtheorem{defi}{\sc Definition}[section]
\newtheorem{theo}[defi]{\sc Theorem}
\newtheorem{prop}[defi]{\sc Proposition}
\newtheorem{nota}[defi]{\sc Notation}
\newtheorem{coro}[defi]{\sc Corollary}
\newtheorem{rema}[defi]{\sc Remark}
\newtheorem{lemm}[defi]{\sc Lemma}
\newtheorem{exem}[defi]{\sc Example}

\newenvironment{demo}[1][]{\noindent {\it Proof.} \protect\nopagebreak
       \rm #1}{\protect\nopagebreak $\square $ \par\stl}

\begin{document}

{ \Huge On meromorphic functions defined by a differential
system of order $1$, II}

\medskip
\begin{center}
{\large Tristan Torrelli}\footnote{Laboratoire 
J.A. Dieudonn\'e, UMR
du CNRS 6621, Universit\'e de Nice Sophia-Antipolis,
Parc Valrose, 06108 Nice Cedex 2, France. {\it E-mail:} 
tristan$\_$torrelli@yahoo.fr


}
\end{center}

\bigskip

\ 

\stl

\noindent{\sc Abstract}. Given a nonzero germ $h$ of 
holomorphic function on $({\bf C}^n,0)$, we study the 
condition: \textsl{``the ideal $\mbox{Ann}_{\cal D}\,1/h$ is generated 
by operators of order $1$''}. When $h$ defines a generic
arrangement of hypersurfaces with an isolated singularity,
we show that it is verified if and only if $h$ is weighted
homogeneous and $-1$ is the only integral root of
its Bernstein-Sato polynomial. When $h$ is a
product, we give a process to test this last condition.
Finally, we study some other related conditions.

\bigskip

\ 

\stl

\section{Introduction}

Let $h\in {\mathcal O}={\mathbf C}\{x_1,\ldots,x_n\}$ be a nonzero germ
of holomorphic function such that $h(0)=0$. We denote by ${\mathcal O}[1/h]$
the ring ${\mathcal O}$ localized by the powers of~$h$. Let
${\mathcal D}={\mathcal O}\langle \partial_1,\ldots, \partial_n\rangle$ 
be the ring of linear differential operators with holomorphic coefficients
and  $F_\bullet{\cal D}$ its filtration by order.
In \cite{TO3}, we study the following condition on $h$:

 \begin{description}
  \item[\ \ A($1/h$)\,:] The left ideal 
 $\mbox{Ann}_{\cal D}\,1/h\subset {\cal D}$ of operators
annihilating $1/h $ is generated by operators of order one.
\end{description}

This property is very natural when one considers 
 sections of  ${\mathcal O}[1/h]/{\mathcal O}$
with an algebraic viewpoint, see \cite{T1}. On the
other hand, it seems to be linked to the topological property 
{\bf LCT($h$)}: \textsl{the de Rham complex 
$\Omega^\bullet[1/h]$ of meromorphic forms with poles 
along $h=0$ is quasi-isomorphic to its subcomplex of 
logarithmic forms}. In particular, {\bf LCT($h$)} implies 
{\bf A($1/h$)} for free germs \cite{Dual} (in the sense of 
K. Saito \cite{Saitlog}). The study of this condition {\bf LCT($h$)} 
was initiated in \cite{CMN} by F.J. Castro Jim\'enez, 
D. Mond and L. Narv\'aez Macarro (see \cite{TO4} for a survey). 
In this paper, we pursue the study of the condition
{\bf A($1/h$)}, and more precisely when $h$ is a reducible germ. 
Our motivation is to deepen the link between  {\bf LCT($h$)} and 
{\bf A($1/h$)}.

\stl

Let us recall that this last condition is closely linked to the
following ones:

 \begin{description}
   \item[\ \ H($h$)\,:] The germ $h$ belongs to the ideal of its partial derivatives.
   \item[\ \ B($h$)\,:] $-1$ is the smallest integral root of the Bernstein
    polynomial of $h$.
  \item[\ \ A($h$)\,:] The ideal $\mbox{Ann}_{\cal D}\,h^s$ is generated
  by operators of order one.
\end{description}

Indeed, condition {\bf H($h$)} seems to be necessary in order to have
{\bf A($1/h$)}, see \cite{TO4}. Moreover, condition {\bf A($1/h$)} 
always implies {\bf B($h$)} (\cite{TO3}, Proposition 1.3).
This last condition has the following algebraic meaning:
\textsl{the ${\cal D}$-module  ${\cal O}[1/h]$ is generated by $1/h$} (see below).
On the other hand, one can easily check that:
\begin{equation} \label{recipABH}
 \mbox{\textsl{If conditions {\bf H($h$)}, {\bf B($h$)} and
{\bf A($h$)} are verified, then so is} 
{\bf A($1/h$)}.}
\end{equation}

Our first part is devoted to condition {\bf B($h$)}. For testing this
condition, it seems natural to avoid the full determination
of the Bernstein polynomial of $h$, denoted by $b(h^s,s)$. Here
we give such a trick when $h$ is not irreducible, using  
Bernstein polynomials associated with sections of holonomic 
${\mathcal D}$-modules.

\stl

Given a nonzero germ $f\in {\mathcal O}$ and an  element $m\in {\mathcal M}$ 
of a holonomic ${\mathcal D}$-module without $f$-torsion, we recall that
 there exists a functional equation:
\begin{equation} \label{eqBernmf}
   b(s)mf^s=P(s)\cdot mf^{s+1}
\end{equation}
in $({\cal D}m)\otimes{\cal O}[1/f,s]f^s$, where $P(s)\in{\cal D}[s]={\cal D}
\otimes {\bf C}[s]$ and $b(s)\in{\bf C}[s]$ are nonzero \cite{K2}. The
{\em Bernstein polynomial} of $f$ associated with $m$, denoted by $b(mf^s,s)$,
is the monic polynomial $b(s)\in{\bf C}[s]$ of smallest degree which verifies such
an equation. When $f$ is not a unit and 
$m\in f^{r-1}{\mathcal M}- f^{r}{\mathcal M}$ with
$r\in {\bf N}^*$, it is easy to check that $-r$ is a root 
of $b(mf^s,s)$. Thus we consider the following condition:
 \begin{description}
   \item[\ \ B($m,f$)\,:] $-1$ is the smallest integral root of $b(mf^s,s)$
\end{description}
for $m\in {\mathcal M}-f{\mathcal M}$; this extends our previous notation
 when $m=1\in{\mathcal O}={\mathcal M}$. 
By generalizing a well known result due to M. Kashiwara, 
this condition means: \textsl{the ${\cal D}$-module $({\cal D}m)[1/f]$ is 
generated by $m/f$} (see Proposition \ref{propmhs}). Hence we get:
\begin{prop} \label{Is-1}
  Let $h_1,h_2\in {\mathcal O}$ be two nonzero germs without common
  factor and such that $h_1(0)=h_2(0)=0$.

\noindent (i) We have: {\bf B($h_1h_2$)}\,$\Rightarrow$\,{\bf B($1/h_1,h_2$)}\,$\Rightarrow$\,{\bf B(}$\mbox{\em\.1}/h_1,h_2${\bf)} where
  $\mbox{\em\.1}/h_1\in {\mathcal O}[1/h_1]/{\mathcal O}$.

\noindent (ii) If {\bf B($h_1$)} is verified, then
{\bf B($h_1h_2$)}\,$\Leftrightarrow$\,{\bf B($1/h_1,h_2$)}.

\noindent (iii) If {\bf B($h_2$)} is verified, then
 {\bf B($1/h_1,h_2$)}\,$\Leftrightarrow$\,{\bf B(}$\mbox{\em\.1}/h_1,h_2${\bf)}.
\end{prop}

Of course, the equivalence  in (ii) just means: $({\mathcal O}[1/h_1])[1/h_2]
={\mathcal O}[1/h_1h_2]$. 
Let us insist on the  condition
{\bf B(}$\mbox{\.1}/h_1,h_2${\bf)}. Indeed, the polynomial
$b((\mbox{\.1}/{h_1})h_2^s,s)$ may be considered as a Bernstein
polynomial of the function $h_2$ in restriction to the
hypersurface $(X_1,0)\subset({\bf C}^n,0)$ defined by $h_1$, see \cite{T1}. 
In particular, $b((\mbox{\.1}/{h_1})h_2^s,s)$ 
coincides with the (classical) Bernstein Sato polynomial of 
$h_2|_{X_1}:(X_1,0)\rightarrow ({\bf C},0)$ if $h_1$
   defines a  smooth germ $(X_1,0)$ (Corollary \ref{Kliss}); thus
this trick is very relevant when $h$ has smooth components.  
As an application, we prove that {\bf B($h$)}
is true when $h$ defines a hyperplane arrangement (Proposition 
\ref{arrhyp}), by using the classical principle of `Deletion-Restriction'. 
 This result was first obtained by A. Leykin \cite{Walt}, and more
recently by M. Saito \cite{MS2}.

\stl

What about the condition {\bf A($1/h$)} when $h=h_1\cdot h_2$ is a 
product  with $h_1(0)=h_2(0)=0$ and $h_1,h_2$ have no common factor ?
It is also natural to consider the ideal $\mbox{Ann}_{\cal D}\,(1/h_1)h_2^s$ and
the Bernstein polynomial $b((1/h_1)h_2^s,s)$. Indeed {\bf B($1/h_1,h_2$)} is a weaker
condition than {\bf B($h_1h_2$)} (Proposition \ref{Is-1}) and we have an analogue
of  (\ref{recipABH}). Of course, it is difficult to verify if 
$\mbox{Ann}_{\cal D}\,(1/h_1)h_2^s$ 
is - or not - generated by operators of order one. Meanwhile,
this may be done under strong assumptions on the components 
of $h$,  by using the characteristic variety of 
${\cal D}(1/h_1)h^s_2$ which may be explicited in terms of the
one of ${\cal D}(1/h_1)$ \cite{Gin}.
Let us give a definition.

\begin{defi}
   A reduced germ $h\in{\cal O}$ defines a {\em generic arrangement of
hypersurfaces with an isolated singularity} if it is a product
$\prod_{i=1}^p h_i$, $p\geq 2$, of germs $h_i$ which defines
an isolated singularity, and such that, for any index $2\leq k\leq \min(p,n)$,
the morphism $(h_{i_1},\ldots,h_{i_k}):({\bf C}^n,0)\rightarrow 
({\bf C}^k,0)$ defines a complete intersection with an isolated singularity at
the origin.
\end{defi}

In the second part, we give a full characterization of {\bf A($1/h$)} for such a 
type of germ.

\begin{theo} \label{carAnnGen}
   Let $h=\prod_{i=1}^ph_i\in{\cal O}$, $p\geq 2$, define a generic
arrangement of hypersurfaces with an isolated singularity. 
Then the ideal
$\mbox{\em Ann}_{\cal D}\,1/h$ is generated by operators of 
order one if and only if the following conditions are verified:
 \begin{enumerate}
    \item the germ $h$ is weighted homogeneous; 
   \item $-1$ is the only integral root of the Bernstein 
 polynomial of $h$. 
 \end{enumerate}
\end{theo}

We recall that a nonzero germ $h$ is {\em weighted homogeneous}
of weight $d\in{\bf Q}^+$ for a system $\alpha\in({\bf Q}^{*+})^n$
if there exists a system of coordinates in which $h$ is a linear
combination of monomials $x_1^{\gamma_1}\cdots x_n^{\gamma_n}$
with $\sum_{i=1}^n\alpha_i\gamma_i=d$.

This result generalizes the case of a hypersurface with an isolated
singularity \cite{T1}. Moreover, the condition {\bf B($h$)}
is also explicit when $p=2$, $h$ weighted homogeneous (Corollary \ref{corpdeux}), 
and the trick above for testing  {\bf B($h$)} may be generalized for $p\geq 3$ 
(Proposition \ref{Bpourplusque2}). On the other hand, 
these conditions on the components of $h$ are strong and 
they are not verified in general. To illustrate this limitation, we end 
this part by studying the condition {\bf A($1/h$)} for $h=(x_1-x_2x_3)g$ 
when $g\in{\bf C}[x_1,x_2]$ is a weighted homogeneous polynomial.

\begin{prop} \label{freex}
   Let $g\in{\bf C}[x_1,x_2]$ be a weighted homogeneous reduced
polynomial of multiplicity greater or equal to $3$. Let
$h\in{\bf C}[x_1,x_2,x_3]$ be the polynomial $(x_1-x_2x_3)g$.

\noindent (i) If $g$ is not homogeneous, then the condition {\bf A($1/h$)} 
does not hold for $h$.

\noindent (ii) If $g$ is homogeneous of degree $3$, then {\bf A($1/h$)} 
holds for $h$.
\end{prop}

Here {\bf H($h$)} are {\bf B($h$)} are verified (see Lemma \ref{courbeplusliss})
whereas {\bf A($h$)} fails. We mention that this family of surfaces was 
intensively studied by the Sevilian group in order to understand 
the condition {\bf LCT($h$)} \cite{3}, \cite{CN}, \cite{CUE}, \cite{CU}, 
\cite{CU1}. 

\stl

In the last part, we give some results on conditions
closely linked to {\bf A($1/h$)}. First, we show how the
Sebastiani-Thom process allows to construct germs
$h$ which verify the condition {\bf A($h$)}. Then, we
do some remarks on a natural generalization of
condition {\bf A($1/h$)}. We end this note with
some remarks on the holonomy of a particular  
$\cal D$-module which appears in the study 
of  {\bf LCT($h$)}.

\stl

\noindent{\bf Aknowledgements.} This research has 
been supported by a Marie Curie Fellowship of the European 
Community (programme FP5, contract HPMD-CT-2001-00097). 
The author is very grateful to the Departamento de \'Algebra, 
Geometr\'{\i}a y Topolog\'{\i}a (Universidad de Valladolid)
 for hospitality during the fellowship, and to the Departamento 
de \'Algebra (University of Sevilla) for hospitality in February 2004
and March 2005.

\section{The condition {\bf B($h$)} for reducible germs}

\subsection{Preliminaries}
 
In this paragraph, we recall some results about Bernstein polynomials of
a germ $f\in{\cal O}$ associated with a section $m$ of a holonomic 
${\cal D}$-module ${\mathcal M}$  without $f$-torsion. As they 
appear in \cite{Ttez} (unpublished), we recall some proofs for the 
convenience of the reader.

\begin{lemm} \label{lembase1}  
  Let $f\in{\cal O}$ be a nonzero germ such that $f(0)=0$. Let 
  $m$  be a germ of holonomic ${\cal D}$-module ${\mathcal M}$
   without $f$-torsion.  Let $P(s)\in{\mathcal D}[s]$ be a differential operator
  such that  $P(j)m f^j\in{\mathcal M}[1/f]$ is zero for a infinite sequence of
  integers  $j\in{\bf Z}$. Then $P(s)$ belongs to the annihilator in
  ${\mathcal D}[s]$ of  $m f^s\in{\mathcal M}[1/f,s]f^s$, denoted by 
  $\mbox{{\em Ann}}_{{\mathcal D}[s]}\,mf^s$.
\end{lemm}

\begin{demo}
  We have the following identity:
\begin{equation}
   \label{lembasic}
  P(s)mf^s=(\sum_{i=0} ^d m_i s^i) f^{s-N}
\end{equation}
 in ${\mathcal M}[1/f,s]f^s$,  where $m_i\in{\mathcal M}$ and $N\in{\bf N}$ 
denotes  the order of $P$. By assumption, there exists some integers 
$j_0<\cdots< j_d$ such that 
   $\sum_{i=0} ^d (j_k)^i m_i=0$ in ${\mathcal M}$ for $0\leq k\leq d$. 
Since the Gram matrix of the integers $j_0,\ldots, j_d$ is inversible, the
previous identities imply that  $m_i=0$ for $0\leq i\leq d$. 
We conclude with (\ref{lembasic}).
\end{demo}

\begin{lemm} \label{lembase2}
    Let $f\in{\cal O}$ be a nonzero germ such that $f(0)=0$.
Let $m\in{\cal M}$ be a nonzero section of a holonomic ${\cal D}$-module 
without $f$-torsion.
  
\noindent (i) If $g\in{\mathcal O}$ is such that  $g\cdot m=0$, then
  $b(mf^s,s)$ coincides with $b(m(f+g)^s,s)$.

\noindent (ii)  If $m\in {\cal M}-f{\cal M}$, then $(s+1)$
divides $b(mf^s,s)$.

\noindent  (iii) For all $p\in{\bf N}^*$, $b(mf^{ps},s)$ divides the
 $\prod_{i=0}^{p-1}b(mf^s,ps+i)$, and the polynomial
 $l.c.m(b(mf^s,ps),\ldots, b(m f^s,ps+p-1))$ divides
  $b(mf^{ps},s)$. In particular, these polynomials have the
  same roots.
\end{lemm}

\begin{demo}
  In order to prove the first point, we just have to check that the
polynomial  $b(m(f+g)^s,s)$ is a multiple of  $b(mf^s,s)$
 for any  $g\in\mbox{Ann}_{\mathcal O}\,m$, and to apply this fact
 with $\tilde{f}=f+g$, $\tilde{g}=-g$. Let $P(s)\in{\mathcal D}[s]$ be a differential
 operator which realizes the Bernstein polynomial of  $m(f+g)^s$. 
 In particular,  $R(s)=b(m(f+g)^s,s)-P(s)f$ belongs to
  $\mbox{Ann}_{{\mathcal  D}[s]}\,m(f+g)^s$. As $(f+g)^j\cdot m=f^j\cdot m$ 
for all $j\in{\bf N}$, the  operator  $R(s)$ annihilates $mf^s$  by
  Lemma \ref{lembase1}. Thus  the polynomial $b(mf^s,s)$ divides
  $b(m(f+g)^s,s)$.

  Now, we prove (ii). Let $R\in {\cal D}$ be the remainder in the division of
$P(s)$ by $(s+1)$ in a nontrivial identity (\ref{eqBernmf}). 
Thus $R\cdot mf^{s+1}=(R\cdot m)f^{s+1}+(s+1)af^s$
where $a\in{\cal M}[1/f,s]$. From (\ref{eqBernmf}), we 
get $b(-1)m=fR(m)$. Hence $b(-1)=0$ since $m\not\in f{\cal M}$.

  The last point is an easy exercice. 
\end{demo}

\begin{prop}
  Let $X\subset {\bf C}^n$ be an analytic subvariety of
codimension $p$ passing through the origin. Let $i:X\hookrightarrow
{\bf C}^n$ denote the inclusion and let $h_1,\ldots, h_p\in {\cal O}$
be local equations of $i(X)$. Let $f\in{\cal O}$ be a germ such that
$f\circ i$ is not constant and let ${\cal M}'$ be a holonomic 
${\cal D}_{X,0}$-module without $(f\circ i)$-torsion. 

If 
$m\in {\cal M}'$ is nonzero,  then 
$b(m(f\circ i)^s,s)$ coincides with the polynomial 
$b(i_+(m)f^s,s)$ where $i_+(m)\in {\cal M}'\otimes
({\cal O}[1/h_1\cdots h_p]/\sum_{i=1}^p{\cal O}[1/h_1\cdots \check{h}_i 
\cdots h_p])$ denotes the element $\mbox{\em\.1}/h_1\cdots h_p$. 
\end{prop}

\begin{demo}
   Up to a change of coordinates, we can  assume that $h_i=x_i$, $1\leq i\leq p$.
Then the remainder $\tilde{f}\in{\bf C}\{x_{p+1},\ldots,x_n\}$ in the
division of $f$ by $x_1,\ldots,x_p$ defines the germ $f\circ i$. Thus we
have $b(i_+(m)f^s,s)=b(i_+(m)\tilde{f}^s,s)$ by using Lemma \ref{lembase2}.
Let us prove that $b(i_+(m)\tilde{f}^s,s)$ coincides with 
$b(m\tilde{f}^s,s)$. Firstly, it easy to check that a functional equation for 
$b(m\tilde{f}^s,s)$ induces an equation for $b(i_+(m)\tilde{f}^s,s)$; thus
$b(i_+(m)\tilde{f}^s,s)$ divides $b(m\tilde{f}^s,s)$. On the other hand, we
consider the following equation:
\begin{equation} \label{equabernlem}
  b(i_+(m)\tilde{f}^s,s)i_+(m)\tilde{f}^s=P\cdot i_+(m)\tilde{f}^{s+1}
\end{equation}
where $P\in{\cal D}[s]$. It may be written $P=\sum_{i=1}^pQ_ix_i+R$
where $Q_i\in {\cal D}[s]$ and the coefficients of $R\in {\cal D}[s]$  
do not depend on $x_1,\ldots, x_p$; in particular, we can change $P$ by 
$R$ in (\ref{equabernlem}). Let $\tilde{R}\in {\cal D}_{X,0}[s]=
{\bf C}\{x_{p+1},\ldots,x_n\}\langle\partial_{p+1},\ldots, 
\partial_n \rangle[s]$ denote the constant term of $R$ as an operator 
in $\partial_1,\ldots, \partial_p$ with coefficients
in ${\cal D}_{X,0}[s]$. Obviously we can change $R$ by $\tilde{R}$ in
(\ref{equabernlem}). As the annihilator of $i_+(m)\tilde{f}^s$ in
${\cal D}_{X,0}[s]$ coincides with the one of $m\tilde{f}^s$, we
deduce that $b(i_+(m)\tilde{f}^s,s)$ is a multiple of $b(m\tilde{f}^s,s)$.
This completes the proof.
\end{demo}

\begin{coro}    \label{Kliss}
   Let   $h_1,h_2\in {\mathcal O}$ be two nonzero germs without common
  factor and such that $h_1(0)=h_2(0)=0$. Assume that $h_1$
   defines a smooth germ $(X_1,0) \subset ({\bf C}^n,0)$. Then 
 $b((\mbox{\em\.1}/{h_1})h_2^s,s)$ coincides with the
   (classical) Bernstein Sato polynomial of
   $h_2|_{X_1}:(X_1,0)\rightarrow ({\bf C},0)$.
\end{coro}

\begin{prop} \label{propmhs}
  Let $f\in{\mathcal O}$ be a nonzero germ such that $f(0)=0$.
   Let  $m$ be a section of a holonomic ${\mathcal D}$-module
    without $f$-torsion, and $\ell\in{\bf N}^*$. The following conditions
   are equivalent:
 \begin{enumerate}
   \item \label{ppetite} The smallest integral root of $b(mf^s,s)$ is
   strictly greater than $-\ell-1$.
   \item \label{engloc} The ${\mathcal D}$-module $({\mathcal D}m)[1/f]$
    is generated by $mf^{-\ell}$.
   \item \label{eval} The following morphism is an isomorphism:
     \begin{eqnarray*}
        \frac{{\cal D}[s]mf^s}{(s+\ell){\cal D}[s]mf^s}
             &\longrightarrow & ({\cal D}m)[1/f] \\
             \stackrel{.}{P(s)mf^s} &\mapsto & P(-\ell)\cdot mf^{-\ell}\ .
     \end{eqnarray*}
\end{enumerate}
\end{prop}

This is a direct generalization of a well known result due to
M. Kashiwara and J.E. Bj\"ork for $m=1\in {\mathcal O}={\mathcal M}$
(see \cite{K1} Proposition 6.2, \cite{Bj} Propositions 6.1.18, 6.3.15 \& 6.3.16).

\subsection{Is $-1$ the only integral root  of $b(h^s,s)$ ?}

\label{CondB}

First of all, let us prove  Proposition  \ref{Is-1}.

\stl

\noindent{\it Proof of Proposition \ref{Is-1}}. 
 Assume that condition {\bf B($h_1h_2$)} is verified. From Proposition
\ref{propmhs}, this means ${\cal D}1/h_1h_2={\cal O}[1/h_1h_2]$. 
In particular, we have $({\cal D}1/h_1)[1/h_2]\subset {\cal D}1/h_1h_2$; 
thus, by using Proposition \ref{propmhs} with $m=1/h_1$, condition 
{\bf B($1/h_1,h_2$)} is verified.
The second relation in (i) is clear since a functional equation
realizing $b((1/h_1)h^s_2,s)$ induces a functional equation
for $b((\mbox{\.1}/{h_1})h_2^s,s)$.

  The second point is clear, since it just means $({\cal O}[1/h_1])[1/h_2]=
{\cal O}[1/h_1h_2]$ (using three times Proposition \ref{propmhs}).
Now,  given $P\in {\cal D}$ and $\ell\in{\bf N}$, let us prove that
$(P\cdot 1/h_1)\otimes 1/h^{\ell}_2$ belongs to ${\cal D}1/h_1h_2$
when {\bf B(}$\mbox{\.1}/h_1,h_2${\bf)} and {\bf B($h_2$)} are verified. From 
Proposition \ref{propmhs}, there exists an operator $Q\in{\cal D}$ 
such that
$(P\cdot\mbox{\.1}/h_1)\otimes 1/h^{\ell}_2=Q\cdot \mbox{\.1}/h_1\otimes 1/h_2$ 
in $({\cal O}[1/h_1]/{\cal O})[1/h_2]$. Hence we have
$(P\cdot 1/h_1)\otimes 1/h^{\ell}_2=Q\cdot 1/h_1h_2 +a/h^N_2$, where
$a\in{\cal O}$ and $N\in{\bf N}^*$. As condition {\bf B($h_2$)} is verified, there
exists $R\in{\cal D}$ such that $R\cdot 1/h_2=a/h_2^N.$ Thus we get
$(P\cdot 1/h_1)\otimes 1/h^{\ell}_2=(Q+Rh_1)\cdot 1/h_1h_2$.
In consequence, the condition {\bf B($1/h_1,h_2$)} is also verified. $\square$

\stl

The following examples show that there is no other relation
between {\bf B($h_1h_2$)}, {\bf B($1/h_1,h_2$)}, {\bf B(}$\mbox{\.1}/h_1,h_2${\bf)}
and {\bf B($h_1$)}, {\bf B($h_2$)}.

\begin{exem}{\em
(i) If $h_1=x_1$ and $h_2=x_1+x_2x_3+x_4x_5$, then
$b(h_1^s,s)=b(h_2^s,s)=s+1$ but
$b((\mbox{\.1}/h_1)h_2^s,s)=b((x_2x_3+x_4x_5)^s,s)=(s+1)(s+2)$ by
using Corollary \ref{Kliss}.

(ii) If $h_1=x_1x_2+x_3x_4$ and $h_2=x_1x_2+x_3x_5$, then
 $b(h_1^s,s)=b(h_2^s,s)=(s+1)(s+2)$, but
 $b((h_1h_2)^s,s)$ is equal to $(s+1)^4(s+3/2)^2$ by using Macaulay 2 
\cite{Mac}, \cite{MacD}.
 Moreover, if $h_3=x_1$, then condition {\bf B($h_1h_3$)} is also true, since
   $b((h_1h_3)^s,s)=(s+1)^3(s+3/2)$ using Macaulay 2. 
Hence condition  {\bf B($h_1h_2$)} does not depend in general of the conditions   
 {\bf B($h_1$)} and {\bf B($h_2$)}.

(iii) Assume that  $h_1=x_1$ and $h_2=x_1^2+x_2^4 +x_3^4$. Then
  $b(h_1^s,s)=s+1$ and condition {\bf B$($}$\mbox{\.1}/h_1,h_2${\bf $)$} 
is true, since $b((\mbox{\.1}/h_1)h_2^s,s)
=b((x_2^4+x_3^4)^s,s)$ by Corollary \ref{Kliss}. But a direct computation 
using \cite{Cras1} shows that condition {\bf B($1/h_1,h_2$)}  is false.

(iv) Assume that $h_1=x_1x_2x_3+x_4x_5$ and $h_2=x_1$. Then
 $ b((1/h_1)h_2^s,s)=b((\mbox{\.1}/h_1)h_2^s,s)=b((x_4x_5)^s,s)=(s+1)^2$, 
using  \cite{18} Proposition 2.9 and  \cite{Cras1} Proposition 1. On the other
 hand, $(s+1)(s+2)$ divides $b((h_1h_2)^s,s)$ and $b(h_1^s,s)$, by the 
 semi-continuity of the Bernstein polynomial (since when $u$ is
a unit, we have $b((u(ux_2x_3+x_4x_5))^s,s)=(s+1)(s+2)$). 
Thus  {\bf B($1/h_1,h_2$)}  does not imply {\bf B($h_1h_2$)} in general.}
\end{exem}

As an application of  Proposition \ref{Is-1}, we obtain a new proof of
 the following result.

\begin{prop} [\cite{Walt}, \cite{MS2}] \label{arrhyp}
   Let $h\in{\bf C}[x_1,\ldots, x_n]$ be the product of nonzero linear forms
(distinct or not). Then the Bernstein polynomial of $h$ has only $-1$ 
as integral root.
\end{prop}

\begin{demo}
  Let $h$ be the product $l_1^{p_1}\cdots l_r^{p_r}$ where
$r,p_1,\ldots,p_r\in{\bf N}^*$ are positive integers, and
$l_i\in({\bf C}^n)^\star$ are distinct. We prove the result
by induction on $r$. If $r=1$, this is a direct
consequence of the following identity:
$$\frac{1}{p^p}\left(\frac{\partial}{\partial x}\right)^p\cdot (x^p)^{s+1}=
(s+\frac{1}{p})(s+\frac{2}{p})\cdots
(s+\frac{p-1}{p})(s+1)(x^p)^s $$
for $p\in {\bf N}^*$. Now, we assume that the assertion is true for any germ
as above with at most $N\geq 1$ distinct irreducible components.
Let $h$ be such a germ with $r=N$. Let $l\in({\bf C}^n)^\star$ be a nonzero
form which is not a factor of $h$, and $p\in {\bf N}^*$. In particular, $-1$
is the only integral root of the Bernstein polynomial of $l$, $l^p$ and $h$.
Let us remark that the assertion for $h\cdot l$ implies the assertion for 
$h\cdot l^p$. Indeed, using Lemma \ref{lembase2}, it is easy to check that 
{\bf B($1/h,l$)} implies {\bf B($1/h,l^p$)}. We conclude with the help 
of Proposition \ref{Is-1}, (ii).

   In order to prove {\bf B($h\cdot l$)}, we just have to check that
$-1$ is the only integral root of $b((\mbox{\.1}/l)h^s,s)$ 
(Proposition \ref{Is-1}, (iii)).
But this is true by induction on $N$ since this last polynomial 
coincides with the Bernstein polynomial of $h|_{\{l=0\}}$ 
(Corollary \ref{Kliss}). This completes the proof.
\end{demo}

When $h$ has more than two components, the following
result provides a generalized criterion for the condition {\bf B($h$)}.

\begin{prop} \label{Bpourplusque2}
   Let $h_1,\ldots, h_p\in{\cal O}$ be nonzero germs without
common factor, and such that $h_1(0)=\cdots = h_p(0)=0$.

(i) Assume that $2\leq p\leq n$ and that $(h_1,\ldots,h_p)$
defines a complete intersection. If  
{\bf B($h_1\cdots \check{h_j}\cdots h_p$)}, $1\leq j\leq p$, are
verified, then {\bf B($\delta,h_1 $)} implies {\bf B($h_1\cdots h_p$)}
where $\delta=\mbox{\em\.1}/h_2\cdots h_p\in {\cal O}[1/h_2\cdots h_p]/
\sum_{i=2}^p{\cal O}[1/h_2\cdots\check{h_i}\cdots h_p]$.

(ii) Assume that $p=n$ and $(h_1,\ldots,h_n)$
defines the origin. If the conditions 
{\bf B($h_1\cdots \check{h_j}\cdots h_n$)}, $1\leq j\leq n$, are
verified, then so is {\bf B($h_1\cdots h_n$)}.

(iii) Assume that $p\geq n+1$. If the conditions 
{\bf B($h_{i_1}\cdots h_{i_n}$)} are verified for
 $1\leq i_1<\cdots <i_n \leq p$ then so is {\bf B($h_1\cdots h_p$)}. 
\end{prop}

\begin{demo}
We start with the first assertion.  From Proposition \ref{Is-1}, we
just have to prove {\bf B($1/h_2\cdots h_p,h_1$)}
(since {\bf B($h_2\cdots h_p$)} is verified).  Thus, given
 $P\in{\cal D}$ and $\ell\in {\bf N}$, let us
prove that $(P\cdot 1/h_2\cdots h_p)\otimes 1/h_1^\ell$ 
belongs to ${\cal D}1/h_1\cdots h_p$. Using condition
{\bf B($\delta,h_1$)}, we have
$$(P\cdot \frac{1}{h_2\cdots h_p})\otimes \frac{1}{h_1^\ell}=
R\cdot \frac{1}{h_1\cdots h_p}+
\sum_{2\leq i\leq p}\frac{q_i}{h_1^{\ell_{i,1}}\cdots
\check{h_i}^{\ell_{i,i}}\cdots h_p^{\ell_{i,p}}}$$
with $q_i\in{\cal O}$ and $\ell_{i,j}\in{\bf N}$. We conclude
by using that  ${\cal O}[1/h_1\cdots \check{h_i}\cdots h_p]$
is generated by $1/h_1\cdots \check{h_i}\cdots h_p$ 
for $2\leq i\leq p$ by assumption.

In order to prove (ii), we have to check that {\bf B($\delta,h_1$)} is
verified when $p=n$. Firstly, we notice that the 
 $\cal D$-module ${\cal O}[1/h_2\cdots h_p]/
\sum_{i=2}^p{\cal O}[1/h_2\cdots\check{h_i}\cdots h_p]$ 
is generated by $\delta$
(using condition {\bf B($h_2\cdots h_p$)}). Thus 
${\cal N}=({\cal D} \delta)[1/h_1]/{\cal D}\delta$ is
isomorphic to the module of local algebraic
cohomology with support in the origin; in particular,
any nonzero section generates $\cal N$. We deduce
easily that $({\cal D}\delta)[1/h_1]$ is generated by
$\delta \otimes 1/h_1$. From Proposition \ref{propmhs},
the condition {\bf B($\delta,h_1$)} is verified.

\stl

The last point is a direct consequence of the following
fact, proved by A. Leykin \cite{Walt}, Remark 5.2:
\textsl{if the condition {\bf B($h_{i_1}\cdots h_{i_{k-1}}$)}
is verified for $1\leq i_1<\cdots <i_{k-1}\leq k$ with  $k\geq n+1$,
then so is {\bf B($h_1\cdots h_k$)}}.
\end{demo}

\begin{exem} {\em
Let  $n=3$, $p\geq 3$ and $h_i=a_{i,1}x^2_1+a_{i,2}x^3_2+a_{i,3}x^4_3$
where the vector $a_i=(a_{i,1},a_{i,2},a_{i,3})$ belongs to ${\bf C}^3$ and the rank
of $(a_{i_1},a_{i_2},a_{i_3})$ is maximal for $1\leq i_1<i_2<i_3\leq p$.
Thus the polynomial $h=h_1\cdots h_p$ defines a generic arrangement
of hypersurfaces with an isolated singularity. By using the closed formulas 
for $b(h^s_i,s)$ and $b((\mbox{\.1}/h_i)h_j^s,s)$, $1\leq i\not=j\leq p$, (see 
\cite{Ya}, \cite{Cras1}), it is easy to check that the  conditions {\bf B($h_i$)}
and {\bf B(}$\mbox{\.1}/h_i,h_j${\bf)} are verified; thus so is {\bf B($h$)}.}
\end{exem}

\section{The condition {\bf A($1/h$)} for a generic arrangement of
hypersurfaces with an isolated singularity} 

In this part, we characterize the condition {\bf A($1/h$)} when 
$h\in{\cal O}$ defines a generic arrangement of hypersurfaces with an 
isolated singularity. Then we study this condition for
a particular family of free germs (\S \ref{freefam}).

\subsection{A convenient annihilator}

This paragraph is devoted to the determination of
an annihilator which will allow us to characterize {\bf A($1/h$)}.

\begin{nota} \label{delta}
   Let $h=(h_1,\ldots,h_r):{\bf C}^n\rightarrow{\bf C}^r$, $1\leq r< n$,
be an analytic morphism. For any $K=
(k_1,\ldots,k_{r+1})\in{\bf N}^{r+1}$ where $1\leq k_1, \ldots,k_{r+1}\leq n$
and $k_i\not=k_j$ for $i\not=j$, let $\Delta^h_K\in{\cal D}$ denote the vector
field: $$\sum_{i=1}^{r+1}(-1)^im_{K(i)}(h)\partial_{k_i}=
 \sum_{i=1}^{r+1}(-1)^i\partial_{k_i} m_{K(i)}(h)$$ where
$K(i)=(k_1,\ldots, \check{k_i},\ldots,k_{r+1})\in{\bf N}^r$ and
$m_{K(i)}(h)$ is the determinant of the $r\times r$ matrix obtained
from the Jacobian matrix of $h$ by deleting the $k$-th columns with
$k\not\in\{k_1,\ldots, \check{k_i},\ldots,k_{r+1} \}$.
\end{nota}

\begin{prop} \label{uncalcgenannu}
Assume that $n\geq 3$.  
Let $h=\prod_{i=1}^p h_i\in{\cal O}$, $p\geq 2$, define a generic
arrangement of hypersurfaces with an isolated singularity, and 
let $\tilde{h}$ be the product $\prod_{i=2}^p h_i$. Then the ideal 
$\mbox{\em Ann}_{\cal D}\,(1/\tilde{h})h_1^s$ is generated by the
operators:
$$\Delta_K^{h_{i_1},\ldots,h_{i_r}}\prod_{i\not=i_1,\ldots,i_r}h_i $$
with $1\leq r\leq \min(n-1,p)$ and $1=i_1<\cdots <i_r\leq p$. 
\end{prop}

\begin{demo}
   Let $I\subset{\cal D}$ be the left ideal generated by the given
operators, and let ${\cal I}\subset {\cal O}[\xi_1,\ldots,\xi_n]$
denote the ideal generated by their principal symbols.
 We will just prove that $\mbox{Ann}_{\cal D}\,(1/\tilde{h})h_1^s
\subset I$, since the reverse inclusion is obvious. Let us study 
$\mbox{char}_{\cal D}\, {\cal D}(1/\tilde{h})h_1^s\subset T^*{\bf C}^n$ 
the characteristic variety of  ${\cal D}(1/\tilde{h})h_1^s$. 
Given an analytic subspace $X\subset {\bf C}^n$, we denote by
$W_{h_1|X}$ the closure in $ T^*{\bf C}^n$ of  the set 
$\{(x,\xi+\lambda dh_1(x))\,|\,\lambda\in{\bf C},\, (x,\xi)\in T^*_X{\bf C}^n \}$.

\stl

\noindent {\it Assertion 1. The characteristic variety of ${\cal D}(1/\tilde{h})h_1^s$ is the
union of the subspaces $W_{h_1}$ and $W_{h_1|X_{i_1,\ldots,i_r}}$, 
$2\leq i_1<\cdots <i_r\leq p$, $1\leq r\leq \min(n-1,p)$, where
$X_{i_1,\ldots,i_r}\subset {\bf C}^n$ is the complete intersection defined by
$h_{i_1},\ldots, h_{i_r}$. }

\stl

\begin{demo}
    Under our assumption, $(\tilde{h}^{-1}(0),x)$ is a germ of a normal crossing 
hypersurface for
any $x\in \tilde{h}^{-1}(0)/\{0\}$ close enough to the origin. In particular, 
${\cal D}1/\tilde{h}$ coincides with ${\cal O}[1/h_{i_1}\cdots h_{i_r}]$
on a neighborhood of such a point, where  
$\{i_1,\ldots,i_r\}=\{i\,|\, h_i(x)=0,\ 2\leq i\leq p\}$.
Hence, the components of the characteristic variety of ${\cal D}1/\tilde{h}$
which are not supported by $h_1=0$ are $T^*_{{\bf C}^n}{{\bf C}^n}$ and
the conormal spaces $T^*_{X_{i_1,\ldots, i_r}}{\bf C}^n$, with 
$2\leq i_1<\cdots <i_r\leq p$ and 
  $1\leq r\leq \min(n-1,p)$. The assertion
follows from a result of V. Ginzburg (\cite{Gin} Proposition 2.14.4). 
\end{demo}

We recall that the relative conormal space\footnote{See \S \ref{SebThom}} 
$W_{h_1}\subset T^*{{\bf C}^n}$ is defined by the polynomials 
$\sigma(\Delta^{h_1}_{k_1,k_2})=h'_{1,x_{k_2}}\xi_{k_1}-h'_{1,x_{k_1}}\xi_{k_2}$,
 $1\leq k_1<k_2\leq n $ (see \cite{Ya} for example). One can also determine
explicitly the defining ideal  of the spaces $W_{h_1|X_{i_1,\ldots, i_r}}$.

\stl

\noindent{\it Assertion 2}\ (\cite{Cras1}). {\it The conormal space 
$W_{h_1|X_{i_1,\ldots, i_r}}$ is defined by $h_{i_1},\ldots, h_{i_r}$ and
  by the principal symbol of the vector fields $\Delta_K^{h_{i_1},\ldots,h_{i_r}}$ 
(when $r<n-1$), where $K=(k_1,\ldots, k_{r+2})\in{\bf N}^{r+2}$
with $1\leq k_1<\cdots <k_{r+2} \leq n$.}

\stl

Now we can determine the equations of 
$\mbox{char}_{\cal D}\, {\cal D}(1/\tilde{h})h_1^s$.

\stl

\noindent{\it Assertion 3. The defining ideal of 
$\mbox{\em char}_{\cal D}\, {\cal D}(1/\tilde{h})h_1^s$ 
is included in $\cal I$.}

\stl

\begin{demo}
  Let $A\in{\cal O}[\xi]={\cal O}[\xi_1,\ldots, \xi_n]$ 
be a polynomial which is zero on the characteristic
variety of  ${\cal D}(1/\tilde{h})h_1^s$. We will prove the result 
when $p\geq n$ - the case $p\leq n-1$ is analogous. 

 Using the inclusion
$W_{h_1|X_{i_1,\ldots,i_{n-1}}}\subset
\mbox{char}_{\mathcal D}\, {\mathcal D}(1/ \tilde{h})h_1^s$ 
and Assertion 2, we have:
$A\in (h_{i_1}, \ldots, h_{i_{n-1}}){\mathcal O}[\xi]$ for
$2\leq i_1<\cdots <i_{n-1}\leq p$. By an easy induction on $p\geq n$, 
one can check that:
$$\bigcap_{2\leq i_1<\cdots <i_{n-1}\leq p}(h_{i_1},\ldots,h_{i_{n-1}}){\mathcal O}
=\sum_{2\leq i_1<\cdots <i_{n-2}\leq p} [\prod_{i\not=1,i_1,\ldots, {i_{n-2}}}h_i ]{\mathcal O}$$
using that every sequence $ (h_{i_1},\ldots,h_{i_{n}})$
is regular. Thus $A$ may be written as a sum 
$\sum_{2\leq i_1<\cdots <i_{n-2}\leq p}A^{(0)}_{i_1,\ldots,i_{n-2}}
(\prod_{i\not=1,i_1,\ldots,i_{n-2}} h_i)$
for some $A^{(0)}_{i_1,\ldots,i_{n-2}}\in {\mathcal O}[\xi]$.

Now let us fix $i_1<\cdots<i_{n-2}$ a family of index as above.  
From  the inclusion $W_{h_1|X_{i_1,\ldots, i_{n-2}}}\subset
\mbox{char}_{\mathcal D}\, {\mathcal D}(1/ \tilde{h})h_1^s$ and Assertion 2,
  $A$ belongs to the ideal
${\cal I}_{1,i_1,\ldots,i_{n-2}}= (h_{i_1}, \ldots, h_{i_{n-2}}){\cal O}[\xi]
+ \sum_K \sigma(\Delta_K^{h_1,h_{i_1},\ldots,h_{i_{n-2}}}) {\mathcal O}[\xi]$.
On the other hand, let us remark that $h_i$ is 
${\cal O}[\xi]/ {\cal I}_{1,i_1,\ldots,i_{n-2}}$-regular for $i\not=1,i_1,\ldots, i_{n-2}$
 [by the principal ideal theorem, using that  ${\cal I}_{1,i_1,\ldots,i_{n-2}}$ defines
the irreducible space $W_{h_1|X_{1,i_1,\ldots,i_{n-2}}}$  of pure dimension $n+1$].
Thus we have
$A^{(0)}_{i_1,\ldots,i_{n-2}}\in  {\cal I}_{1,i_1,\ldots,i_{n-2}}$, and $A$ may 
be written:
$A=U+\sum_{2\leq i_1<\cdots <i_{n-3}\leq p}A^{(1)}_{i_1,\ldots,i_{n-3}}
(\prod_{i\not=1,i_1,\ldots,i_{n-3}} h_i)$
where $A^{(1)}_{i_1,\ldots,i_{n-3}}\in {\mathcal O}[\xi]$  and 
$U\in{\mathcal I}$. Up to a division by ${\mathcal I}$, we can assume that $U=0$.
After iterating this process with  $W_{h_1|X_{i_1,\ldots,i_r}}$,
$1\leq r\leq n-2$, we deduce that $A-A^{(n-2)} \tilde{h}$ 
belongs to ${\mathcal I}$.
Hence, using that $W_{h_1}\subset
 \mbox{char}_{\mathcal D}\, {\mathcal D}(1/ \tilde{h})h_1^s$, we have: 
$A^{(n-2)}\in \sum_{1\leq k_1<k_2\leq n} 
   \sigma(\Delta_{k_1,k_2}^{h_1}) {\mathcal O}[\xi]$. In particular, 
$A^{(n-2)}\tilde {h}$ belongs to ${\mathcal I}$, and we conclude 
that $A\in {\mathcal I}$.
\end{demo}

  Now let us prove the proposition. Let 
$P\in \mbox{Ann}_{\mathcal D}\,(1/ \tilde{h})h_1^s$ be a
nonzero operator of order $d$. In particular, $\sigma(P)$ is zero on
$\mbox{char}_{\mathcal D}\, {\mathcal D}(1/\tilde {h})h_1^s$, and by 
Assertion 3: $\sigma(P)\in{\cal I}$. In other words, there exists
$Q\in I$ such that $\sigma(Q)=\sigma(P)$. Thus, the operator 
$P-Q\in \mbox{Ann}_{\mathcal D}\,(1/ \tilde{h})h_1^s \cap F_{d-1}{\cal D}$
belongs to $I$, and so does $P$ (by induction on 
the order of operators). 
\end{demo}

\begin{rema}{\em
We are not able to determine $\mbox{Ann}_{\cal D}\,h^s$
when $h$ defines a generic arrangement of hypersurfaces
with an isolated singularity. In particular, we do not know
if the condition {\bf A($h$)} (or {\bf W($h$)}) 
is - or not - verified (see \S \ref{SebThom}).}
\end{rema}

Given a germ $h\in{\cal O}$ such that $h(0)=0$, let us
denote by $\mbox{Der}(-\log h)$ the coherent
$\cal O$-module of logarithmic derivations relative to $h$, 
that is, vector fields  which preserve $h{\cal O}$ (see \cite{KS}).

\begin{coro}
    Let $h=\prod^p_{i=1}h_i\in{\cal O}$, $p\geq 2$, define
 a generic arrangement of hypersurfaces with an isolated
 singularity. Assume that $n\geq 3$ and that $h$ is a
 weighted homogeneous polynomial. Then $\mbox{\em Der}(-\log h)$
 is generated by the Euler vector field $\chi$ such that $\chi(h)=h$
and the vector fields
$$\left[\prod_{i\not=i_1,\ldots, i_r}h_i\right]\cdot 
\Delta_K^{h_{1_1},\ldots,h_{i_r}} $$
where $1\leq r\leq \min(n-1,p)$ and $1=i_1<\cdots <i_r\leq p$.
\end{coro}

\begin{demo}
   We denote by $\tilde{h}\in{\cal O}$ the product $h_2\cdots h_p$.
Let $v$ be a logarithmic vector field; in particular, $v(h)=a h$. As
$h=h_1\tilde{h}$, it is easy to check that $v(h_1)=a_1h_1$ and
$v(\tilde{h})=\tilde{a}h_1$ for $a_1,\tilde{a}\in{\cal O}$ such that
$a_1+\tilde{a}=a$. In particular, 
$v\cdot (1/\tilde{h})h_1^s=(a_1s-\tilde{a})(1/\tilde{h})h_1^s$.
Thus $v+\tilde{a}-a_1\chi$ belongs to $\mbox{Ann}_{\cal D}\,(1/\tilde{h})h^s_1$,
and by using the proof of the  previous result, we have:
$$v=-\tilde{a}+a_1\chi+\sum_{r=1}^{\min(n-1,p)}\sum_{1\leq i_1<\ldots< i_r\leq p} 
\lambda_{i_1,\ldots, i_r}\Delta_K^{i_1,\ldots,i_r}\cdot\prod_{i\not=i_1,\ldots, i_r}h_i\ $$
where $\lambda_{i_1,\ldots, i_r}\in{\cal O}$ for $1\leq i_1<\ldots< i_r\leq p$. 
As $v$ is a vector field, we get $v=a_1\chi+\sum_r\sum
\lambda_{i_1,\ldots,i_r}[\prod_{i\not=i_1\cdots i_r}h_i]\Delta_K^{i_1,\ldots,i_r}$ 
and the assertion follows.
\end{demo}

\subsection{The expected characterization}

The proof of Theorem \ref{carAnnGen} is an easy consequence
of the following result

\begin{prop}
   Let $h=\prod_{i=1}^ph_i\in{\cal O}$, $p\geq 2$, define a generic
arrangement of hypersurfaces with an isolated singularity. 
Assume that $n\geq 3$ and that the origin is a critical point
of $h_1$. Let $\tilde{h}$ denote the product 
$\prod_{i=2}^p h_i$. Then the ideal
$\mbox{\em Ann}_{\cal D}\,1/h$ is generated 
by operators of  order one if and only if the following conditions 
are verified:
 \begin{enumerate}
    \item the germ is weighted homogeneous; 
   \item $-1$ is the smallest integral root of the Bernstein 
 polynomial $b((1/\tilde{h})h_1^s,s)$. 
 \end{enumerate}
\end{prop}

\begin{demo}
    We can assume that $h$ does not define a normal crossing divisor. 
Indeed,  the conditions {\bf A($1/h$)}, $1$ and $2$ are obviously 
verified for a normal crossing divisor. In particular, the constant term 
with the coefficient on the right side  of  any operator in 
$\mbox{Ann}_{\cal D}(1/\tilde{h})h^s_1$ is not
a unit (see Proposition \ref{uncalcgenannu}). 

\stl

  Firstly, we prove that conditions $1$ \& $2$ imply {\bf A($1/h$)}.
By an Euclidean division, we have a decomposition
$$\mbox{Ann}_{{\cal D}[s]}\,\frac{1}{\tilde{h}}h^s_1={\cal D}[s](s-\tilde{q}-v)
+{\cal D}[s]\mbox{Ann}_{\cal D}\frac{1}{\tilde{h}}h^s_1$$
where $v$ denotes the Euler vector field such that $v(h_1)=h_1$ and
$v(\tilde{h})=\tilde{q}\tilde{h}$ with $\tilde{q}\in {\bf Q}^{*+}$. 
Moreover, with  the condition $2$, the ideal $\mbox{Ann}_{\cal D}\,
1/(\tilde{h}h_1)$ is obtained by fixing $s=-1$ in a system of generators
of $\mbox{Ann}_{{\cal D}[s]}(1/\tilde{h})h^s_1$ 
(see \cite{T1} Proposition 3.1). From Proposition \ref{uncalcgenannu}, 
the condition {\bf A($1/h$)} is therefore verified.

\stl 

  Now, we prove the reverse. Let us assume that 
$\mbox{Ann}_{{\cal D}}1/h$ is generated by the operators
$Q_1,\ldots,Q_w\in F_1{\cal D}$. From Proposition 1.3 in \cite{TO3},
 {\bf B($h$)} is verified, and so\footnote{In fact, 
 the same proof shows directly that condition {\bf A($1/h$)} 
implies {\bf B($1/\tilde{h},h_1$)}.} is condition $2$ 
by Proposition \ref{Is-1}.  Hence, we just have to check that
$h$ is necessarily weighted homogeneous. Let $q_i$ be the germ
$Q_i(1)\in{\cal O}$ and $Q'_i$ the vector field $Q_i-q_i$. In 
particular, we have $Q'_i(h)=q_ih$ for $1\leq i\leq w$. 
As $h=h_1\tilde{h}$, it is easy to deduce that 
$Q'_i(\tilde{h})=\tilde{q}_i\tilde{h}$ and $Q'_i(h_1)=q_{i,1}h_1$
 where $\tilde{q}_i,q_{i,1}\in {\cal O}$ verify
$$\tilde{q}_i+q_{i,1}=q_i,\ 1\leq i\leq w.$$
On the other hand, we have the following fact:

\stl

\noindent{\it Assertion 1. There exists a differential operator $R$ in 
$\mbox{\em Ann}_{{\cal D}}(1/\tilde{h})h^s_1$ such that 
$R=1+\sum_{i=1}^wA_iq_{i,1}$ with $A_i\in{\cal D}$.}  

\stl

\begin{demo}
  The proof is analogous to the one of \cite{T1} Lemme 3.3. 
From \cite{Gin} p 351 or \cite{Ttez}, there exists 
 a `good' operator $R_0(s)$ of degree $N\geq 1$ in
$\mbox{Ann}_{{\cal D}[s]}(1/\tilde{h})h^s_1$, that is 
$R_0(s)=s^N+\sum_{k=0}^{N-1}s^kP_k$ with $P_k\in F_{N-k}{\cal D}$,
$0\leq k\leq N-1$. By Euclidean division, 
we have $R_0(s)=(s+1)S(s)+R_0(-1)$
where $S(s)$ is monic in $s$ of degree $N-1$ and 
$R_0(-1)\in\mbox{Ann}_{\cal D}1/h$. Thus, there exists 
$A_1,\ldots, A_w\in {\cal D}$ such that $R_0(-1)=\sum_{i=1}^w A_iQ_i$. 
From the relations above, we get
$$(s+1)S(s)\frac{1}{\tilde{h}}h_1^s+(s+1)\sum_{i=1}^wA_iq_{i,1}\frac{1}{\tilde{h}}h^s=0.$$
Hence $R_1(s)=S(s)+\sum_{i=1}^wA_i q_{i,1}$ belongs to
$\mbox{Ann}_{{\cal D}[s]}(1/\tilde{h})h_1^s$. By iteration, we can assume
that $S(s)=1$.
\end{demo}

In particular, at least one of the $q_{i,1}$ is a unit (see the very
 beginning of the proof.)

\stl

\noindent{\it Assertion 2. If $q_{i,1}$ is a unit, then so is $q_i$.}

\stl

\begin{demo}
   As the assertion is clear if $\tilde{q}_i$ is not a unit, 
we can assume that $\tilde{q}_i$ is a unit. Let $\chi_i$ denote
the vector field $q_{i,1}^{-1}Q'_i$; in particular $\chi_i(h_1)=h_1$. As $h_1$
defines an isolated singularity, a famous result due to K. Saito \cite{KS}
asserts that, up to a change of coordinates, $\chi_i$ is an Euler 
vector field $\sum_{k=1}^n\alpha_kx_k\partial_k$ with 
$\alpha_k\in {\bf Q}^{*+}$. Hence, the relation 
$\chi_i(\tilde{h})=q_{i,1}^{-1}\tilde{q}_i\tilde{h}$ implies that 
the constant $(q_{i,1}^{-1}\tilde{q}_i)(0)$ belongs to ${\bf Q}^{*+}$
[consider the initial part of $q_{i,1}^{-1}\tilde{q}_i\tilde{h}$ relative to
$\alpha_1,\ldots, \alpha_n$]. In particular, $q_{i,1}^{-1}\tilde{q}_i+1$
is a unit, and so is $q_i=\tilde{q}_i+q_{i,1}$. 
\end{demo}

We recall that a formal power series $g\in{\bf C}[[x_1,\ldots, x_n]]$
is {\em weakly weighted homogeneous} of type 
$(\beta_0,\beta_1,\ldots, \beta_n)\in {\bf C}^{n+1}$ if for all
monomial $x_1^{\gamma_1}\cdots x_n^{\gamma_n}$ with a nonzero 
coefficient in the power expansion of $g$, we have 
$\beta_1\gamma_1+\cdots+\beta_n{\gamma_n}=\beta_0$. Let us
 pursue the proof. We have proved that there exists an Euler
vector field $\chi_i$ 
such that $q_i^{-1}\chi_i(h)=h$ (in particular, $q_i(0)>0$). From 
\cite{KS}, Corollary 3.3, there exists a formal change of coordinates
$\phi$ such that $h\circ\phi$ is weakly weighted homogeneous
of type $(1, \alpha_1q^{-1}_i(0),\ldots,\alpha_nq^{-1}_i(0))$. As the 
$\alpha_kq^{-1}_i(0)$ are strictly positive, 
$h\circ\phi$ is in fact weighted homogeneous, and according to
a theorem of Artin \cite{Art}, a convergent change of coordinates 
exists. This completes the proof.
\end{demo}

\noindent{\it Proof of Theorem \ref{carAnnGen}.} 
The case $n=2$ is done in \cite{T1}, Theorem 1.2.
We just have to check that the condition $2$ in the
previous statement may be replaced by {\bf B($h$)}.
Indeed, condition {\bf A($1/h$)} always implies
{\bf B($h$)} (\cite{TO3} Proposition 1.3), and
on the other hand, {\bf B($h$)} is stronger than
{\bf B($1/\tilde{h},h_1$)} (Proposition \ref{Is-1}). $\square$

\stl

  Of course, we can use \S\ref{CondB} to test if  condition {\bf B($h$)}
is verified. In the particular case $p=2$ and $h$ weighted homogeneous, 
we obtain the following characterization:

\begin{coro} \label{corpdeux}
  Let $h_1,h_2\in{\bf C}[x_1,\ldots, x_n]$ be two weighted 
homogeneous polynomial of degree $d_1,d_2$ for a system 
$\alpha\in({\bf Q}^{*+})^n$, 
defining hypersurfaces with an isolated singularity at the origin and 
 without common components.
Let ${\cal K}\subset{\cal O}$ be the ideal generated by the 
maximal minors of the Jacobien matrix of $(h_1,h_2)$. Then 
the annihilator of $1/h_1h_2$ is generated by operators of order 
$1$ if and only if  for $j=1$ or $2$, 
 there is no weighted homogeneous element in 
${\cal O}/h_j{\cal O}+{\cal K}$
whose weight belongs to the set $\{d_j\times k -\sum_{i=1}^n \alpha_i\ ;\ 
k\in{\bf N}\ \&\ k\geq 2 \}$. 
\end{coro}

This relies on the existence of  closed formulas for
$b((1/\tilde{h})h^s_1,s)$ under these assumptions 
\cite{Cras1}.


\subsection{About a family of free germs} \label{freefam}

In this part, we prove Proposition \ref{freex}. As the two parts are quite
distinct, we will prove them successively.

\begin{lemm}\label{courbeplusliss}
Let $g\in {\bf C}\{x_1,x_2\}$ be a nonzero reduced germ of plane curve
 such that $g(0)=0$. Then $-1$ is the only integral root of the Bernstein
polynomial of $(x_1-x_2x_3)g(x_1,x_2)$.
\end{lemm}

\begin{demo}
  As $g$ is a reduced germ of plane curve, 
{\bf B($g$)} is verified \cite{20}, \cite{15}. Thus, by using
Proposition \ref{Is-1}, the three conditions 
{\bf B($(x_1-x_2x_3)g(x_1,x_2)$)}, {\bf B($1/x_1-x_2x_3,g$)}
and {\bf B(}$\mbox{\.1}/x_1-x_2x_3,g${\bf)} are equivalent. Let
us prove the last one. From Corollary \ref{Kliss}, we have
$b((\mbox{\.1}/x_1-x_2x_3)g^s,s)=b((g(x_2x_3,x_2))^s,s)$. 
Let us write $g(x_2x_3,x_2)=x_2^\ell\tilde{g}(x_2,x_3)$ where
$\tilde{g}\in {\bf C}\{x_2,x_3\}- x_2{\bf C}\{x_2,x_3\}$ is reduced
and $\ell \in {\bf N}^*$. If $\tilde{g}$ is a unit, then
{\bf B($g(x_2x_3,x_3)$)} is verified and so is 
 {\bf B($(x_1-x_2x_3)g(x_1,x_2)$)}. Now we assume that $\tilde{g}$
is not a unit. As it is reduced, {\bf B($\tilde{g}$)} is verified and
{\bf B($\tilde{g}x_2^\ell$)} is equivalent to {\bf B($1/\tilde{g},x_2^\ell$)}.
Using Lemma \ref{lembase2}, it is easy to check that 
{\bf B($1/\tilde{g},x_2$)} implies  {\bf B($1/\tilde{g},x_2^\ell$)}. 
Thus we just have to prove {\bf B($1/\tilde{g},x_2$)}. As
condition {\bf B$(\tilde{g})$} is verified, the conditions 
{\bf B($1/\tilde{g},x_2$)}, {\bf B($\tilde{g}x_2$)}
and {\bf B(}$\mbox{\.1}/x_2,\tilde{g}${\bf)} are equivalent
 (Proposition \ref{Is-1}). Both of them are verified since  
$b((\mbox{\.1}/x_2)\tilde{g}^s,s)=b((\tilde{g}(0,x_3))^s,s)$
 from Corollary \ref{Kliss}, where $\tilde{g}(0,x_3)=ux_3^N$
with $u\in{\bf C}\{x_3\}$ is a unit. This completes the proof.
\end{demo}

We recall that a nonzero germ $h\in{\cal O}$ defines a
germ of {\it free} divisor if the module of logarithmic
derivations relative to $h$ is ${\cal O}$-free
\cite{Saitlog}. Moreover, such a germ defines a 
{\em Koszul-free} divisor if there exists a basis 
$\{\delta_1,\ldots, \delta_n\}$ of 
$\mbox{Der}(-\log h)$ such that the sequence of
principal symbols 
$(\sigma(\delta_1),\ldots,\sigma(\delta_n))$ is 
$\gr^F {\cal D}$-regular.

\begin{lemm} \label{hfree}
   Let $g\in {\bf C}[x_1,x_2]$ be a weighted homogeneous and 
reduced polynomial whose multiplicity is greater or equal to $3$.
Let $h\in {\bf C}[x_1,x_2,x_3]$ denote  the polynomial 
$(x_1- x_2 x_3)g(x_1,x_2)$.

  (i) The polynomial $h$ defines a free divisor and verifies
the condition \mbox{\bf H($h$)}.
  
  (ii) The polynomial $h$ defines a Koszul-free divisor
 if and only if the weighted homogeneous
 polynomial $g$ is not homogeneous.
\end{lemm}

\begin{demo}
(i)  It is enough to remark that the following vector fields
verify Saito's criterion \cite{Saitlog} for $h$:
\begin{eqnarray*}
  \delta_1 & = & \alpha_1 x_1 \partial_1 + \alpha_2 x_2 \partial_2
+(\alpha_1-\alpha_2)x_3 \partial_3 \\
  \delta_2 & = & g'_{x_2}\partial_1-g'_{x_1}\partial_2+(x_3u-v)\partial_3\\
  \delta_3 & = & (x_1-x_2x_3)\partial_3
\end{eqnarray*}
where $(\alpha_1,\alpha_2)\in ({\bf Q}^{*+})^2$ is a system of weights for $g$,
and  $u\in {\bf C}[x_1,x_2,x_3]$, $v\in {\bf C}[x_2,x_3]$
are the polynomials of degree in $x_3$ less or equal to $1$
uniquely defined by the relation
$$x_3g'_{x_1}(x_1,x_2)+g'_{x_2}(x_1,x_2)=u(x_1,x_2,x_3)x_1-v(x_2,x_3)x_2$$
(we use that $g'_{x_1},g'_{x_2}\in(x_1,x_2){\bf C}[x_1,x_2]$ under our
assumptions.)

\stl

\noindent (ii) As the sequence $(\sigma(\delta_1),\sigma(\delta_2),\xi_3)$ 
is regular, the
germ $h$ is Koszul-free if and only if the sequence
$(\sigma(\delta_1),\sigma(\delta_2), x_1-x_2 x_3)$ is 
${\cal  O}[\xi]$-regular.
By division by $x_1-x_2 x_3$, this condition may be rewritten:
\textsl{the polynomials
\begin{eqnarray*}
   \Upsilon_1 & = & \alpha_1 x_2 x_3 \xi_1 +\alpha_2x_2
\xi_2+(\alpha_1-\alpha_2)x_3\xi_3 \\
 \Upsilon_2 & = & g'_{x_2}(x_2 x_3, x_2)\xi_1-g'_{x_1}(x_2 x_3,x_2)\xi_2+
(x_3u(x_2 x_3,x_2,x_3)-v(x_2,x_3))\xi_3 
\end{eqnarray*}
have no common factor.} Let us notice that $x_2$ is the only (irreducible) 
common factor of $g'_{x_1}(x_2 x_3,x_2)$ and $g'_{x_2}(x_2 x_3,x_2)$ 
[since $g\in {\bf C}[x_1,x_2]$ defines an isolated singularity.] Thus, when 
$\Upsilon_1$ and $\Upsilon_2$  have a  common factor, this factor is $x_2$  
(up to a multiplicative constant).
As $g$ belongs in $(x_1,x_2)^3{\bf C}[x_1,x_2]$, we have
$g'_{x_1},g'_{x_2}\in(x_1,x_2)^2{\bf C}[x_1,x_2]$; thus
$u,v\in(x_1,x_2){\bf C}[x_1,x_2,x_3]$. In
particular, $x_2$ is a factor of $\Upsilon_2$, and $\Upsilon_1,
\Upsilon_2$ have no common factor if and only if
$\alpha_1\not=\alpha_2$. This completes the proof.
\end{demo}

Of course, for $g=x_1x_2(x_1+x_2)$, $h$ is the example of
F.J. Calder\'on-Moreno in \cite{3} and it is not Koszul-free.

\stl

\noindent{\it Proof of Proposition \ref{freex}, part (i).}
Without loss of generality, we will assume that $\delta_1(h)=h$.
Let us take $\delta'_2=\delta_2-u\cdot\delta_1$ and
$\delta'_3=\delta_3+x_2\delta_1$; in particular,
$\{\delta_1,\delta'_2,\delta'_3\}$ is a basis of
Der(log $h$) such that $\delta'_2(h)=\delta'_3(h)=0$.

From the characterization of condition {\bf A($1/h$)} for
Koszul-free germs (see \cite{TO3} Corollary 1.8), it is enough
to check that condition {\bf A($h$)} fails, that is,
the sequence $(x_1-x_2x_3,\sigma(\delta'_2),\sigma(\delta'_3))$ 
is not regular. As $g$ belongs to $(x_1,x_2)^3{\bf C}[x_1,x_2]$,
we have $\sigma(\delta'_2),\sigma(\delta'_3) \in (x_1,x_2){\cal O}[\xi]$.
By division by $x_1-x_2x_3$, we deduce that the sequence is
not regular. $\square$

\begin{nota} {\em
  Given a homogeneous polynomial $g\in{\mathbf C}[x_1,x_2]-{\bf C}$ of
degree $p\geq 1$, we denote by $\tilde{g}_1, \tilde{g}_2\in
{\mathbf C}[x_1,x_2,x_3]$ the quotient of the division of
$g'_{x_1}, g'_{x_2}$ by $x_1-x_2x_3$. In particular:}
\begin{equation}
   \label{gtilde}
  g'_{x_i}=(x_1-x_2x_3)\tilde{g}_i+x_2^{p-1}g'_{x_i}(x_3,1),\ i\in\{1,\,2\}.
\end{equation}
\end{nota}

\begin{lemm}
  Let $g\in{\mathbf C}[x_1,x_2]$ be a homogeneous reduced polynomial
of degree $p\geq 3$. Then the characteristic variety of
${\cal D}(1/x_1-x_2x_3)g^s$ is defined by the following polynomials:
$(x_1-x_2x_3)\xi_3$, $g'_{x_2}\xi_1-g'_{x_1}\xi_2+px_2^{p-2}
g(x_3,1)\xi_3$, and $[x_2g'_{x_2}(x_3,1)\xi_1-x_2g'_{x_1}(x_3,1)\xi_2+pg(x_3,1)\xi_3]\xi_3.$
\end{lemm}

\begin{demo}
  Using \cite{Gin} Proposition 2.14.4, the characteristic variety
of  the $\cal D$-module  ${\cal D}(1/x_1-x_2x_3)g^s$ is the union of 
the conormal spaces
$W_g$ and $W_{g|x_1=x_2x_3}$. It is easy to check that they are defined
by the ideals $I_1=(\xi_3,g'_{x_2}\xi_1-g'_{x_1}\xi_2){\mathcal O}[\xi]$ and
$I_2=(x_1-x_2x_3,x_2g'_{x_2}(x_3,1)\xi_1-x_2g'_{x_1}(x_3,1)\xi_2+pg(x_3,1)\xi_3)
{\mathcal O}[\xi]$ respectively. Clearly, the ideal $I$ generated by
the given polynomials is contained in $I_1\cap I_2$. Thus we just have to prove
the reverse relation.

\stl

  Let $A,B,C,D\in{\mathcal O}[\xi]$ be such that
$$ A(x_1-x_2x_3)+B(x_2g'_{x_2}(x_3,1)\xi_1-x_2g'_{x_1}(x_3,1)\xi_2+pg(x_3,1)\xi_3)=   
C\xi_3+D(g'_{x_2}\xi_1-g'_{x_1}\xi_2).$$
Using (\ref{gtilde}), we get
$$(A-D(\tilde{g}_2\xi_1-\tilde{g}_1\xi_2))(x_1-x_2x_3)+(pBg(x_3,1)-C)\xi_3 $$
$$+\ (B-Dx_2^{p-2})x_2(g'_{x_2}(x_3,1)\xi_1-g'_{x_1}(x_3,1)\xi_2)=0$$
Since the sequence 
$(x_1-x_2x_3,\xi_3,x_2(g'_{x_2}(x_3,1)\xi_1-g'_{x_1}(x_3,1)\xi_2))$
is ${\mathcal O}[\xi]$-regular, there exist $U,V,W\in{\mathcal O}[\xi]$ such
that
$$\left\{\begin{array}{rcl}
    A-D(\tilde{g}_2\xi_1-\tilde{g}_1\xi_2)&=&U\xi_3+Wx_2(g'_{x_2}(x_3,1)\xi_1-g'_{x_1}(x_3,1)\xi_2)\\
    B-Dx_2^{p-2}&=&-V\xi_3-W(x_1-x_2x_3)
         \end{array}  \right.$$
Thus one can notice that the first part of the first identity belongs to $I$,
that is,  $I$ is the defining ideal of $W_g\cup W_{g|x_1=x_2x_3}$.
\end{demo}

\begin{lemm}
  Let $g\in{\mathbf C}[x_1,x_2]$ be a homogeneous reduced polynomial
of degree $3$. Then the annihilator of $(1/x_1-x_2x_3)g^s$ is generated
by the following differential operators:
$$(x_1-x_2x_3)\partial_{3}-x_2\mbox{\,,\ \  }g'_{x_2}\partial_1-g'_{x_1}\partial_2+3x_2
g(x_3,1)\partial_3+ x_3\tilde{g}_1+\tilde{g}_2 \mbox{\ \ }and$$
$$[x_2g'_{x_2}(x_3,1)\partial_1-x_2g'_{x_1}(x_3,1)\partial_2+3g(x_3,1)\partial_3]
\partial_3+\tilde{g}_2\partial_1 -\tilde{g}_1\partial_2 + 3g'_{x_1}(x_3,1)\partial_3+u'_{x_1}$$

\vspace{2mm}
\noindent where $u=x_3\tilde{g}_1+\tilde{g}_2$.
\end{lemm}

\begin{demo}
  Let us denote by $I\subset {\cal D}$ the ideal generated by
the given operators $S_1$, $S_2$, $S_3$. It is not hard to check
the inclusion $I\subset\mbox{Ann}_{\cal D}\,(1/x_1-x_2x_3)g^s$.
Let us prove that the reverse inclusion by induction on the order
of operators.

\stl

  Let $P\in\mbox{Ann}_{\cal D}\,(1/x_1-x_2x_3)g^s$ be an operator of
order $d$. As $d=0$ implies $P=0$, we can assume $d\geq 1$. Then
$\sigma(P)$ is zero on the characteristic variety of 
 ${\cal D}(1/x_1-x_2x_3)g^s$.
From the previous result, there exists $A_1\in{\mathcal O}[\xi]$ (resp.
$A_2$, $A_3$) zero or homogeneous in $\xi$ of degree $d-1$ (resp. $d-1$, $d-2$)
such that: $\sigma(P)=\sum_{i=1}^3A_i\sigma(S_i)$. If $\tilde{A}_i\in{\cal D}$,
 $1\leq i\leq 3$, are such that $\sigma(\tilde{A}_i)=A_i$ for $1\leq i\leq 3$,
then $P-\sum_{i=1}^3\tilde{A}_iS_i$ belongs to $F_{d-1}{\cal D}$ and annihilates
$(1/x_1-x_2x_3)g^s$. By induction, it belongs to $I$ and so does $P$.
\end{demo}

\noindent {\it Proof of Proposition \ref{freex}, part (ii).} We will prove that
$\mbox{Ann}_{\cal D}\,1/h$ is generated by the operators 
$\tilde{\delta}_1=\delta_1+4$,
$\tilde{\delta}_2=\delta_2+u$, $\tilde{\delta}_3=\delta_3-x_2$ (with
the notations introduced in the proof of Lemma \ref{hfree} with 
$\alpha_1=\alpha_2=1$).
From Lemma \ref{courbeplusliss}, we know that $-1$ is the smallest integral
 root of $b((1/x_1-x_2x_3)g^s,s)$.
Thus we have the decomposition
$\mbox{Ann}_{\cal D}\,1/h={\cal D}\tilde{\delta}_1+\mbox{Ann}_{\cal D}\,(1/x_1-x_2x_3)g^s$,
and the assertion is a direct consequence of the previous result and of the
following relation in ${\cal D}$:
$$[g'_{x_2}(x_3,1)x_2\partial_1-g'_{x_1}(x_3,1)x_2\partial_2+3g(x_3,1)\partial_3 +3g'_{x_1}(x_3,1)](\partial_3\tilde{\delta}_1-\partial_1\tilde{\delta}_3)$$
$$+\ [\partial_2+x_3\partial_1](\partial_3\tilde{\delta}_2 + (\tilde{g}_2\partial_1-\tilde{g}_1\partial_2)\tilde{\delta}_3)\ =\ -2 S_3+\partial_1
\tilde{\delta}_2-(\tilde{g}_2\partial_1-\tilde{g}_1\partial_2+u'_{x_1})\tilde{\delta}_1$$
where $S_3$ is the operator of order $2$ which appears in the given system
of generators of $\mbox{Ann}_{\cal D}\,(1/x_1-x_2x_3)g^s$. $\square$

\section{Some other conditions}

In this part, $h\in{\cal O}$ denotes a nonzero germ such 
that $h(0)=0$. 

\subsection{The condition A($h$) for Sebastiani-Thom germs} 
\label{SebThom}

We recall that the condition {\bf A($h$)} on the ideal
$\mbox{Ann}_{\cal D}\, h^s$ may be considered
almost as a geometric condition. Indeed the
following condition implies {\bf A($h$)}:
\begin{description}
 \item[\ \ W($h$)\,:] The relative conormal space
$W_h$ is defined by linear equations in $\xi$.
\end{description}
since $W_h=\overline{\{(x,\lambda dh)\,|\, \lambda\in{\bf C}\}} 
\subset T^*{\bf C}^n$  is the characteristic variety of 
${\cal D}h^s$ (\cite{K1}). For example, {\bf W($h$)} is
true for  hypersurfaces with an isolated singularity \cite{Ya} and
for locally weighted homogeneous free divisors \cite{CN}.
This condition also means that the kernel of the morphism
of graded $\cal O$-algebras:
\begin{eqnarray*}
 {\cal O}[X_1,\ldots, X_n]& \longrightarrow &
{\cal R}({\cal J}_h) \\
 X_i & \longmapsto &  th'_{x_i} 
\end{eqnarray*}
is generated by homogeneous elements of degree $1$, where
${\cal J}_h$ denotes the Jacobian ideal $(h'_{x_1},\ldots,
h'_{x_n}){\cal O}$ and ${\cal R}({\cal J}_h)$ is the Rees
algebra $\bigoplus_{d\geq 0}{\cal J}_h^dt^d$. Following a
terminology due to W.V. Vasconcelos, one says that ${\cal J}_h$
is \textsl{of linear type} (see \cite{CN} for more details). Finally, let
us give a third condition trapped between {\bf A($h$)} and
{\bf W($h$)}:
\begin{description}
 \item[\ \ G($h$)\,:] The graded ideal $\gr^F\mbox{Ann}_{\cal D}\,h^s$
is generated by homogeneous polynomials in $\xi$ of degree $1$.
\end{description}

\begin{rema}{\em
(i) We do not know if the conditions {\bf A($h$)}, {\bf G($h$)} and {\bf W($h$)}  
are - or not - equivalent.

\noindent (ii) These conditions are not stable by multiplication by a unit.}
\end{rema}

It seems uneasy to find sufficient conditions on $h$ for {\bf A($h$)} or
{\bf W($h$)}. Thus, it is natural to study if the class of germs $h$ which
 verify {\bf A($h$)} or {\bf W($h$)} is - or not - stable by Thom-Sebastiani 
sums. Here we give a positive answer in a particular case.

\begin{prop}
  Let $g\in{\cal O}$ be a nonzero germ such that $g(0)=0$ and
which verifies the condition {\bf W($g$)}. Let $f\in{\bf C}\{z_1,\ldots, z_p\}$ 
be a nonzero germ which defines an isolated singularity at the origin. 
Then $h=g+f$ verifies the condition {\bf W($h$)}. 
\end{prop}

  This is direct consequence of the following result.

\begin{prop}
  Let $g\in{\cal O}$ be a nonzero germ such that $g(0)=0$,
and $\Upsilon_1,\ldots,\Upsilon_w\in{\cal O}[\xi]$ be
homogeneous polynomials which generate the defining
ideal of $W_g$.

Let $f\in{\bf C}\{z_1,\ldots, z_p \}$ be a nonzero germ which defines an
isolated singularity and $\xi_1,\ldots, \xi_n,\eta_1,\ldots,\eta_p$ denote
the conormal coordinates on $T^*{\bf C}^n\times {\bf C}^p$.  Then
 the relative conormal space 
$W_{g+f}\subset T^*{\bf C}^n\times {\bf C}^p$ is defined by the
polynomials $f'_{z_i}\eta_j-f'_{z_j}\eta_i$, $1\leq i<j\leq p$, 
 $g'_{x_k}\eta_i-f'_{z_i}\xi_k$, $1\leq i\leq p$, $1\leq k\leq n $, and 
$\Upsilon_1,\ldots,\Upsilon_w$.
\end{prop}

\begin{demo}
  Let us denote by $E\subset {\bf C}\{z_1,\ldots,z_p\}$ a 
${\bf C}$-vector space of finite dimension isomorphic to 
${\bf C}\{z_1,\ldots,z_p\}/(f'_{z_1},\ldots, f'_{z_p})$ by projection, and
by ${\bf C}\{x,z\}$ the ring ${\bf C}\{x_1,\ldots,x_n,z_1,\ldots,z_p\}$. 
In particular, any germ $p\in{\bf C}\{x,z\} $ may be written in a 
unique way: $p=\tilde{p}+r$ where 
$\tilde{p}\in E\otimes_{\bf C}{\cal O}\subset {\bf C}\{x,z\}$ and
$r\in  (f'_{z_1},\ldots, f'_{z_p}){\bf C}\{x,z\}$. 

We denote by $I_{f+g}\subset  {\bf C}\{x,z\}[\xi,\eta]$ 
the ideal generated by the given operators,
and by $I_g\subset {\bf C}\{x,z\}[\xi,\eta]$ (resp. $I_f$) the
ideal generated by $\Upsilon_1,\ldots, \Upsilon_w$ 
(resp. $f'_{z_i}\eta_j-f'_{z_j}\eta_i$, $1\leq i<j\leq p$). Obviously,
any element of $I_{g+f}$ is zero on $W_{g+f}$. Let us prove the
reverse relation.

\stl

  Let $P\in {\bf C}\{x,z\}[\xi,\eta]$ be a homogeneous polynomial of
degree $N\in{\bf N}^*$ in $(\xi,\eta)$ which is zero on $W_{g+f}$.

\stl

\noindent {\it Assertion 1. There exists 
$\tilde{P}(\xi,\eta)\in {\bf C}\{x,z\}[\xi,\eta]$ such that
$P-\tilde{P}(\xi,\eta)$ belongs to $I_{g+f}$, and it is of the form:
$$\tilde{P}(\xi,\eta)=Q(\eta)+\sum_{|\gamma|\leq N-1} 
\tilde{P}_\gamma(\xi)\eta_1^{\gamma_1}\cdots \eta_p^{\gamma_p}$$
where $\gamma=(\gamma_1,\ldots,\gamma_p)\in{\bf N}^p$,
 $\tilde{P}_\gamma(\xi)\in(E\otimes {\cal O})[\xi]$ are zero 
or homogeneous in $\xi$ of degree $N-|\gamma|$, 
$Q(\eta)\in {\bf C}\{x,z\}[\eta]$ is zero or homogeneous of degree $N$.}

\stl

\begin{demo}
   Let us write: $P=\sum_{|\beta+\gamma|=N}p_{\beta,\gamma}\eta^\gamma
\xi^\beta$ with $p_{\beta,\gamma}\in{\cal O}$. For all $\beta\in{\bf N}^n$,
$|\beta|=N$, the germ $p_{\beta,0}$ may be written in a unique
way $p_{\beta,0}=\tilde{p}_{\beta,0}+r_{\beta,0}$ with 
$\tilde{p}_{\beta,0}\in E\otimes {\cal O}$
and $r_{\beta,0}=\sum_{i=1}^pr_{\beta,0,i}f'_{z_i} $ for some 
$r_{\beta,0,i}\in{\bf C}\{x,z\}$. As $|\beta|\geq 1$, there exists 
an index $k$ such that $\beta_k\not=0$. Thus
$$r_{\beta,0}\xi_1^{\beta_1}\cdots \xi_n^{\beta_n}
-\sum_{i=1}^p r_{\beta,0,i}g'_{x_k}\eta_i\xi_1^{\beta_1}\cdots 
 \xi_k^{\beta_k-1}\cdots \xi_n^{\beta_n}  \in I_{g+f}$$
and we fix $\tilde{P}_0(\xi)=\sum_{|\beta|=N}\tilde{p}_{\beta,0}\xi^\beta$.
By iterating this process for increasing $|\gamma|$, we get a decomposition
$P=Q(\eta)+\sum_{|\gamma|\leq N-1} 
\tilde{P}_\gamma(\xi)\eta^{\gamma}+R$ where $R\in I_{g+f}$.
\end{demo}

\noindent {\it Assertion 2. The polynomials $\tilde{P}_\gamma(\xi)$ belong to $I_g$.}

\stl

\begin{demo}
    We prove it by induction on $\gamma$, using the lexicographical order
on ${\bf N}^p$. 
  As $\tilde{P}(g'_{x_1},\ldots, g'_{x_n},f'_{z_1},\ldots,f'_{z_p})=0$, we have
$\tilde{P}_0( g'_{x_1},\ldots, g'_{x_n})\in(f'_{z_1},\ldots, f'_{z_p}){\bf C}\{x,z\}$.
Thus $\tilde{P}_0(\xi)$ belongs to $I_g$ (since 
$\tilde{P}_0(\xi)\in(E\otimes {\cal O})[\xi]$ and $g\in {\cal O}$). Now,
let us assume that $\tilde{P}_{{\gamma}'}(\xi)\in I_g$ for all 
${\gamma}'<\gamma$, ${\gamma}'\geq 0$ and 
$\tilde{P}_{{\gamma}}(\xi)\not=0$.
Since $\tilde{P}(g'_{x_1},\ldots, g'_{x_n},f'_{z_1},\ldots,f'_{z_p})=0$ and
 $\tilde{P}_{{\gamma}'}(g'_{x_1},\ldots, g'_{x_n})=0$  
for ${\gamma}'<\gamma$, we have:
\begin{eqnarray*}
\tilde{P}_\gamma(g'_{x_1},\ldots,g'_{x_n}){f'}_{z_1}^{\gamma_1}\cdots
 {f'}^{\gamma_p}_{z_p}&\in& ({f'}_{z_1}^{\gamma_1+1}, 
{f'}_{z_1}^{\gamma_1}{f'}_{z_2}^{\gamma_2+1},\ldots, 
{f'}_{z_1}^{\gamma_1}\cdots {f'}_{z_{p-1}}^{\gamma_{p-1}}{f'}_{z_p}^{\gamma_p+1})
{\bf C}\{x,z\}\\
& & \mbox{\ \ \ \ \ \ \ \ \ }+\  Q(f'_{z_1},\ldots, f'_{z_p}){\bf C}\{x,z\} \\
& & \mbox{\ \ } \subset ({f'}^{\gamma_1+1}_{z_1},\ldots,{f'}^{\gamma_p+1}_{z_p}  ){\bf C}\{x,z\}
\end{eqnarray*}
since the degree of $Q(\eta)$ is strictly greater than $|\gamma|$.
From this identity, we deduce that 
$\tilde{P}_\gamma(g'_{x_1},\ldots, g'_{x_n})\in(f'_{z_1},\ldots, f'_{z_p}){\bf C}\{x,z\}$
using that $(f'_{z_1},\ldots, f'_{z_p})$ is a ${\bf C}\{x,z\}$-regular  sequence. Thus
$\tilde{P}_\gamma(\xi)$ belongs to $I_g$ as above.
\end{demo}

In particular, the polynomial $P-Q(\eta)$ belongs to $I_{g+f}$. As
$P$ is zero on $W_{g+f}$, we have $Q(f'_{z_1},\ldots, f'_{z_p})=0$.
Thus $Q(\eta)$ belongs to $I_f$ (since $(f'_{z_1},\ldots, f'_{z_p})$ is
${\bf C}\{x,z\}$-regular). We conclude that $P\in I_{g+f}$, 
and this completes the proof. 
\end{demo}

\begin{rema}{\em
Let us recall that the reduced Bernstein polynomial of the germ $h=g(x)+z^N$ 
has no integral root for $N$ `generic' \cite{15}. In particular, our result
allows to construct some examples of weighted homogeneous polynomials
$h$ which verify condition {\bf A($1/h$)}  [with the help of 
identity (\ref{recipABH}) of the Introduction].}
\end{rema}

\subsection{The condition A$_{\log}$($1/h$)}

Let us recall how the condition {\bf A($1/h$)} appears in the study
of the so-called logaritmic comparison theorem.
If $D$ is a free divisor, F.J. Calder\'on-Moreno and L. Narv\'aez-Macarro
\cite{Dual} have obtained a differential analogue of the condition
{\bf LCT($D$)}; in particular, it implies that the natural $\cal D$-linear morphism
$\varphi_D:{\cal D}_X\otimes_{{\cal V}_0^D}{\cal O}_X(D)\longrightarrow
{\cal O}_X(\star D)$ is an isomorphism. Here
${\cal O}_X(D)$ denotes the ${\cal O}_X$-module of meromorphic functions
with at most a simple pole along $D$, and ${\cal V}_0^D\subset{\cal D}_X$
is the sheaf of ring of logarithmic differential operators, that is,
$P\in{\cal D}_X$ such that $P\cdot (h_D)^k\subset (h_D)^k{\cal O}$ for any 
$k\in {\bf N}$, where $h_D$ is a (local) defining equation of $D$.
Locally, we have ${\cal O}_X(D)={\cal V}_0^D\cdot(1/h_D)$, thus
$\varphi_D$ is given by
\begin{eqnarray*}
{\cal D}/{{\cal D}\mbox{Ann}_{{\cal V}_0^D}\,1/h_D} & 
\longrightarrow & {\cal O}[1/h_D] \\
P &\longmapsto& P\cdot \frac{1}{h_D}
\end{eqnarray*}
where $\mbox{Ann}_{{\cal V}_0^D}\,1/h_D\subset{\cal V}_0^D$ 
is the ideal of logarithmic operators which annihilate
$1/h_D$. From the structure theorem of logarithmic operators
associated with a free divisor \cite{3}, we have 
${\cal V}_0^D={\cal O}_X[\mbox{Der}(-\log\,h_D)]$; hence
the ideal $\mbox{Ann}_{{\cal V}_0^D}\,1/h_D$ is locally generated by
$v_i+a_i$, $1\leq i\leq n$, where $\{v_1,\ldots, v_n\}$ is a
basis of $\mbox{Der}(-\log\,h_D)$ and $a_i\in{\cal O}$ is
defined by $v_i(h_D)=a_ih_D$, $1\leq i\leq n$. In particular,
the injectivity of $\varphi_D$ means that the condition 
 {\bf A($1/h$)} is verified.

\stl

  Let us notice that the following  condition may
also be considered:
 
\begin{description}
   \item[\ \ A$_{\log}$($1/h$)\,:] The ideal $\mbox{Ann}_{\cal D}\,1/h$ is 
generated by logarithmic operators.
\end{description}

In this paragraph, we compare these two conditions. Firstly, it is
easy to see that the condition  {\bf A($1/h$)} always implies {\bf A$_{\log}$($1/h$)}.
On the other hand, we do not know if these conditions
are distinct or not. Meanwhile, we have the following result:

\begin{lemm}
  Let $h\in{\cal O}$ be a nonzero germ such that $h(0)=0$. Assume that
one of the following conditions is verified:
\begin{enumerate}
  \item the ring ${\cal V}_0^D$ coincides with ${\cal O}[\mbox{\em Der}(-\log\,h)]$, 
the $\cal O$-subalgebra of ${\cal D}$ 
   generated by the logarithmic derivations relative to $h$.
  \item the conditions {\bf A($h$)} and {\bf H($h$)} are verified. 
\end{enumerate}
Then the conditions {\bf A($1/h$)} and {\bf A$_{\log}$($1/h$)} are equivalent.
\end{lemm}

\begin{demo}
  Assume that condition $1$ is verified, and let $P\in{\cal V}_0^D\cap\mbox{Ann}_{\cal D}\,1/h$
be a nonzero logarithmic operator annihilating $1/h$. By assumption, it may be written as a
sum $\sum_{|\gamma|\leq d}p_\gamma v_1^{\gamma_1}\cdots v_N^{\gamma_N}$
where $p_\gamma\in{\cal O}$ and $v_1,\ldots,v_N$ is a generating system of 
$\mbox{Der}(-\log\,h)$. As $\mbox{Der}(-\log\,h)$ is stable by Lie
brackets, we have 
$$P=\sum_{|\gamma|\leq d}p_\gamma (v_1+a_1)^{\gamma_1}\cdots (v_N+a_N)^{\gamma_N}
+\underbrace{\sum_{|\gamma|<d}r_\gamma v_1^{\gamma_1}\cdots v_N^{\gamma_N}}_R$$
where $r_\gamma\in{\cal O}$, and $a_i\in{\cal O}$ is defined by $v_i(h)=a_ih$, $1\leq i \leq N$;
in particular, $R$ belongs to ${\cal V}_0^D\cap\mbox{Ann}_{\cal D}\,1/h$. By induction,
we conclude that $P$ belongs to the ideal ${\cal D}(v_1+a_1,\ldots, v_N+a_N)$; thus 
{\bf A$_{\log}$($1/h$)} implies the condition {\bf A($1/h$)}.

\stl

 Now we assume that the conditions {\bf A$_{\log}$($1/h$)}, {\bf A($h$)} and
 {\bf H($h$)} are verified. From Proposition \ref{AnnLogB}, the condition
 {\bf B($h$)} is also verified. Thus so is  {\bf A($1/h$)} (see (\ref{recipABH})
in the Introduction). This completes the proof.
\end{demo}

In particular, these conditions coincides for weighted homogeneous polynomials
which define an isolated singularity.

\begin{rema} {\em
 Some criterions for condition $1$ are given by M. Schulze in \cite{Schul}.}    
\end{rema}

Finally, this condition {\bf A$_{\log}$($1/h$)} always implies
{\bf B($h$)} (as {\bf A($1/h$)} does.) 

\begin{prop} \label{AnnLogB}
  Let $h\in{\cal O}$ be a nonzero germ such that $h(0)=0$. If the ideal
$\mbox{\em Ann}_{\cal D}\,1/h$ is generated by logarithmic operators, then
$-1$ is the only integral root of the Bernstein polynomial of $h$.    
\end{prop}

\begin{demo}
  The proof is analogous to the one of \cite{T1}, Proposition 1.3. We
need the following fact.

\stl

\noindent{\em Assertion. If $Q$ is a logarithmic operator relative to
$h$, then $Q\cdot h^s=q(s)h^s$ with $q(s)\in{\cal O}[s]$.}

\begin{demo}
    We have $Q\cdot h^s=a(s)h^{s-N}$ with
$a(s)=\sum_{i=0}^Na_is^i$, $a_i\in {\cal O}$, and $N$ is the degree of $Q$. 
Thus we just have to prove that 
$a(s)\in h^N{\cal O}[s]$. As $Q$ is logarithmic,
$Q\cdot h^k$ belongs to $h^k{\cal O}$ for $k\geq 1$; in
particular $\sum_{i=0}^N a_i k^i\in h^N{\cal O}$ for $1\leq k\leq N+1$. 
By solving this system, we get $a_i\in h^N{\cal O}$, 
$0\leq i\leq N$, that is, $a(s)\in h^N{\cal O}[s]$. 
\end{demo}

Let $Q_1,\dots, Q_w$ be a generating system of logarithmic operators
which annihilate $1/h$. For $1\leq i\leq w$, we have $Q_i\cdot h^s=q_i(s)h^s$
with $q_i(s)\in {\cal O}[s]$. As $Q_i$ annihilates $1/h$, the polynomial 
$q_i(s)$ belongs to $(s+1){\cal O}[s]$ and we
denote $\tilde{q}_i(s)\in{\cal O}[s]$ the quotient of $q_i(s)$ by $(s+1)$.
Let us suppose that the
Bernstein polynomial of $h$, denoted by $b(s)$, has an integral root
strictly smaller than $-1$. We denote by
$k\leq -2$, the greatest integral root of $b(s)$ verifying this
condition. Using a Bernstein equation which gives $b(s)$, we
get:
$$ b(s)\cdots b(s-k-2)h^s=P(s) h^{s-k-1}$$
where $P(s)\in{\mathcal D}[s]$.  Thus $P(k)$ annihilates $1/h$
 and it may be written  $\sum_{i=1}^w A_i Q_i$
  with  $A_i\in {\mathcal D}$, $1\leq i\leq w$. If $P'(s)\in{\mathcal D}[s]$ 
is the quotient of  $P(s)$ by $s-k$,  the previous equation becomes:
 $$\underbrace{b(s)\cdots b(s-k-2)}_{c(s)}h^s=
(s-k)\left[P'(s)+\sum_{i=1}^wA_i\tilde{q}_i\right]h^{-k-2}
\cdot h^{s+1}$$
where $-k-2\geq 0$ and the multiplicity of  $k$ in $c(s)$ is the same in
$b(s)$. Hence, by division by $(s-k)$, we get
a Bernstein functional equation   
such that the polynomial in the left member is not a multiple of 
$b(s)$. But this is not possible, because $b(s)$ is the 
Bernstein polynomial of $h$. Hence we have the result. 
\end{demo}

\subsection{The condition {\bf M($h$)}}

  Let $h\in{\cal O}$ be a nonzero germ such that $h(0)=0$.
In this paragraph, we study the following condition
 \begin{description}
   \item[\ \ M($h$)\,:] The $\cal D$-module 
 $\widetilde{\cal M}_h={\cal D}/\tilde{I}_h$ is holonomic
\end{description}
where $\tilde{I}_h\subset {\cal D}$ is the left ideal generated
by the operators of order $1$ which annihilate $1/h$. 
This condition only depends on the ideal $h{\cal O}$ 
(since the right multiplication by a unit $u\in{\cal O}$ 
induces an isomorphism of $\cal D$-modules
from $\widetilde{\cal M}_h$ to $\widetilde{\cal M}_{uh}$). 

\stl

Let us recall that this condition and this `logarithmic' 
$\cal D$-module - introduced by F.J Castro-Jim\'enez and
J.M. Ucha in \cite{7} -  
are very natural in this topic. Indeed, the condition {\bf A($1/h$)} 
always implies {\bf M($h$)}, since {\bf A($1/h$)}
means that the morphism $\widetilde{\cal M}_h\rightarrow {\cal O}[1/h]$
defined by $P \mapsto P\cdot 1/h$ is an isomorphism. Moreover, 
the condition {\bf LCT($D$)} needs locally {\bf M($h_D$)} for a free
divisor $D$  (see the beginning of the previous paragraph).

\stl

  Here, we link the condition {\bf M($h$)} with some other
conditions introduced in this topic (see \S \ref{SebThom}). Firstly, let us 
consider the following one:
 \begin{description}
   \item[\ \ L($h$)\,:] The ideal in ${\cal O}_{T^*{\bf C}^n}$ 
generated by $\pi^{-1}\mbox{Der}(-\mbox{log\,}h) $
defines an analytic space of  (pure) dimension $n$ 
\end{description}
\noindent where  $\pi$ denotes the canonical map 
$T^*{\bf C}^n\rightarrow {\bf C}^n$. In K. Saito's language,
one says that the irreducible components of the {\em logarithmic
characteristic variety} are holonomic; moreover, this is
equivalent to the local finiteness of the logarithmic
stratification associated with $h$ (see \cite{Saitlog}, \S 3).
 For a free germ, this is exactly the notion of Koszul-free germ 
(see \cite{Saitlog}; \cite{Bru}, Proposition 6.3; \cite{CN}, Corollary 1.9).

\begin{prop}
   Let $h\in {\mathcal O}$ be a nonzero germ such that $h(0)=0$.

\noindent (i) The condition {\em \bf L($h$)} implies {\em\bf M($h$)}.

\noindent (ii) The condition {\em\bf A($h$)} implies {\em\bf M($h$)}.

\noindent (iii) The condition {\em\bf G($h$)} implies {\em\bf L($h$)}.

\noindent (iv) If $h$ defines a locally weighted homogeneous
divisor, then the condition  {\em\bf L($h$)} is verified.
\end{prop}

\begin{demo}
  The first point is clear since 
 $\pi^{-1}\mbox{Der}(-\mbox{log\,}h)\subset \mbox{gr\,}\tilde{I}_h$. 
 Let us prove (ii). By assumption, 
the ideal $J=\mbox{Ann}_{\mathcal D}\, h^s$ is
included $\tilde{I}$. On the other hand, it is obvious that the
operators $h\partial_i+h'_{x_i}$, 
$1\leq i\leq n$, belong to  $\tilde{I}$. 
Hence, we have the following inclusion:
$\mbox{gr}^F J+(h\xi_1,\ldots, h\xi_n){\mathcal O}[\xi]\subset
 \mbox{gr}^F \tilde{I}$. We notice that
 $$\mbox{gr}^FJ+(h\xi_1,\ldots, h\xi_n){\mathcal O}[\xi]=
(\mbox{gr}^FJ,h){\mathcal O}[\xi] \cap (\xi_1,\ldots, \xi_n){\mathcal O}[\xi]  $$
 since $\mbox{gr}^FJ\subset (\xi_1,\ldots, \xi_n){\mathcal O}[\xi]$.
 Thus the characteristic variety of $\widetilde{{\mathcal M}}_h$
 is included in 
$V( \mbox{gr}^FJ,h) \cup V(\xi_1,\ldots, \xi_n)\subset T^*{\bf C}^n$.
 Let us recall that the characteristic variety of  ${\mathcal D}h^s$ is the
the closure $W_h\subset T^*{\bf C}^n$ of
the set $\{(x,\lambda dh(x))\, |\, \lambda \in{\bf C} \}$ \cite{K1}; 
in particular, $W_h$ is irreducible of pure dimension  $n+1$.
From the principal ideal theorem, $W_h\cap \{h=0\}=V(\mbox{gr}^FJ,h)$
has a pure dimension $n$. Hence $\widetilde{{\mathcal M}}_h $
is holonomic.

\stl
 
  The proof of (iii) is the very same, since the ideal
generated by the principal  symbol of the elements in 
$\mbox{Der}(-\log\,h)$ contains
$\mbox{gr}^F J+(h\xi_1,\ldots , h\xi_n){\cal O}[\xi]$.

\stl

Let us prove (iv), by induction on dimension. Let $D\subset {\bf C}^n$ 
denote the hypersurface defined by $h$, and let $L$ be the associated
 logarithmic characteristic variety. If $n=2$, then
{\bf W($h$)} is verified and so is {\bf L($h$)} by (iii). Now, we assume
that $n\geq 3$. From Proposition 2.4 in \cite{CMN}, there exists a
neighborhood $U$ of the origin such that, for each point $w\in U\cap D$,
$w\not=0$, the germ of pair $({\bf C}^n,D,w)$ is isomorphic to a product
$({\bf C}^{n-1}\times{\bf C},D'\times {\bf C},(0,0))$ where $D'$ is a 
locally weighted homogeneous divisor of dimension $n-2$. Up to
this identification, $\mbox{Der}(-\mbox{log\,}h)_w$ is generated
by the elements in  $\mbox{Der}(-\mbox{log\,}h_{D'})$ and
$\partial/\partial z$, where $z$ is the last coordinate on
${\bf C}^{n-1}\times {\bf C}$; in particular, the induction 
hypothesis applied to $D'$ implies the
result for ${\bf C}\times D'$. Hence, the dimension of 
$L\cap \pi^{-1}(U-\{0\})= L-T^*_{\{0\}}{\bf C}^n$ is $n$.
Let $C\subset L$ be an irreducible component of $L$.
If $\pi(C)=\{0\}$, then $C$ coincides with $T^*_{\{0\}}{\bf C}^n$
since $\dim C$ is at most equal to $n$ 
(see \cite{Bru}, Proposition 1.14 (i)). Now, if $\pi(C)$ is not
the origin, then 
$\dim C=\dim (C-T^*_{\{0\}}{\bf C}^n)= 
\dim (L-T^*_{\{0\}}{\bf C}^n)=n$. We conclude that 
$L$ has dimension $n$. 
\end{demo}

  We recall that  K. Saito proved that the condition {\bf L($h$)} 
is verified for any hyperplane arrangements \cite{Saitlog}, Example 3.14. 
The point (iv) may be considered as a  generalization of this  fact.
On the other hand, it generalizes also the fact that
locally weighted homogeneous free divisors are
Koszul-free \cite{LQH}  (of course, our proof is similar).

\stl

The following diagram summarizes the previous relations:

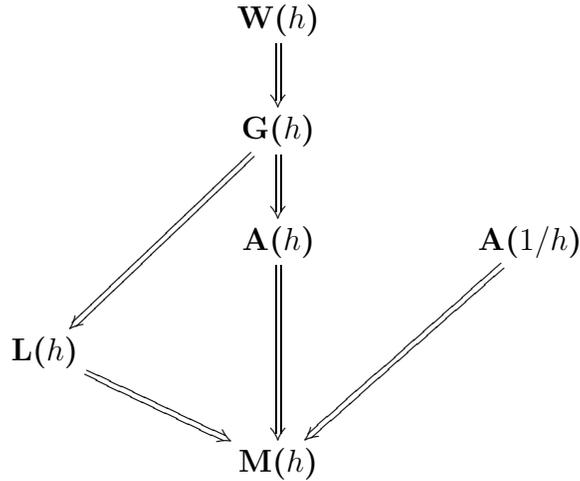
\stl
\begin{displaymath}
\xymatrix{
 &  &\mbox{\bf W($h$)} \ar@{=>}[d] & & \\
 & & \mbox{\bf G($h$)} \ar@{=>}[ddll] \ar@{=>}[d]& & \\
 & & \mbox{\bf A($h$)} \ar@{=>}[dd] & & 
\mbox{\bf A($1/h$)}\ar@{=>}[ddll] \\
     \mbox{\bf L($h$)}\ar@{=>}[drr] &  & & & \\
&  & \mbox{\bf M($h$)} & &}
\end{displaymath}
\stl

Let us notice that the reverse relations are false. Firstly, 
if $h$ is the germ $(x_1-x_2x_3)(x_1x^2_2+x_1^2x_2)$ then 
{\bf L($h$)} and {\bf A($h$)} are not verified but {\bf A($1/h$)}
holds \cite{Saitlog}, \cite{CCMN}, \cite{CN}, \cite{CUE}, \cite{TO3}. On the other hand, if $h=(x_1-x_2x_3)(x_1^3+x^4_2)$
then it defines a Koszul-free germ (see Lemma \ref{hfree} for instance);
in particular, {\bf L($h$)} is verified where as {\bf A($h$)} and 
{\bf A($1/h$)} fail (see the proof of Proposition \ref{freex}, (i)). 
Finally, L. Narv\'aez-Macarro and F.J Calder\'on-Moreno prove in 
\cite{Dual} that the free divisor defined by
$h=(x_1-x_2x_3)(x_1^5+x_2^4+x_1^4x_2)$ is not of Spencer 
type\footnote{This is a necessary condition on a free divisor $D$
for verifying {\bf LCT($D$)}, see \cite{Dual}.}.
In fact, the condition {\bf M($h$)} is no more verified, since all
elements of a system of generators of $\tilde{I}$ belongs 
to ${\cal D}(x_1,x_2)$, see \cite{Dual} \S 5.

\end{document}